\documentclass[12pt, reqno]{amsart}
\usepackage{amsmath}
\usepackage{amssymb}

\numberwithin{equation}{subsection}

\begin{document}

{\theoremstyle{plain}
    \newtheorem{theorem}{\bf Theorem}[subsection]
    \newtheorem{proposition}[theorem]{\bf Proposition}
    \newtheorem{claim}[theorem]{\bf Claim}
    \newtheorem{lemma}[theorem]{\bf Lemma}
    \newtheorem{corollary}[theorem]{\bf Corollary}
}
{\theoremstyle{remark}
    \newtheorem{remark}[theorem]{\bf Remark}
    \newtheorem{example}[theorem]{\bf Example}
}
{\theoremstyle{definition}
    \newtheorem{defn}[theorem]{\bf Definition}
    \newtheorem{question}[theorem]{\bf Question}
}


\newcommand{\lreg}[1]{{\operatorname{fiber-reg}_{#1}}}
\newcommand{\la}[1]{a_{[#1]}}
\newcommand{\rreg}[1]{{{\operatorname{res-reg}_{#1}}}}
\newcommand{\ppi}[1]{\pi^{[#1]}}
\newcommand{\X}[1]{{X^{[#1]}}}
\newcommand{\mS}[1]{S^{[#1]}}
\newcommand{\lF}[1]{F^{[#1]}}
\newcommand{\kset}[1]{[#1]}
\renewcommand{\c}[1]{c^{[#1]}}
\newcommand{\lr}[1]{r_{#1}}
\newcommand{\supp}[1]{{\operatorname{Supp}^{[#1]}}}

\def\to{\longrightarrow}
\def\height{\operatorname{ht}}
\def\reg{\operatorname{reg}}
\def\Reg{\operatorname{\bf reg}_B}
\def\Rreg{\operatorname{\bf res-reg}}
\def\Lreg{\operatorname{\bf fiber-reg}}
\def\depth{\operatorname{depth}}
\def\Hom{\operatorname{Hom}}
\def\proj{\operatorname{Proj}}
\def\grade{\operatorname{grade}}
\def\spec{\operatorname{Spec}}
\def\Tor{\operatorname{Tor}}
\def\Ext{\operatorname{Ext}}
\def\suppo{\operatorname{Supp}}
\def\Ann{\operatorname{Ann}}
\def\La{{\frak a}_{\mathbb F}}

\def\mm{{\frak m}}
\def\pp{{\frak p}}
\def\qq{{\frak q}}
\def\O{{\mathcal O}}
\def\I{{\mathcal I}}
\def\M{{\frak M}}
\def\A{{\mathcal A}}
\def\B{B}
\def\L{{\mathcal L}}
\def\NN{{\mathbb N}}
\def\N{{\mathcal N}}
\def\PP{{\mathbb P}}
\def\ZZ{{\mathbb Z}}
\def\GG{{\mathbb G}}
\def\T{{\mathcal T}}
\def\F{{\mathcal F}}
\def\G{{\mathcal G}}
\def\H{{\mathcal H}}
\def\U{{\frak U}}
\def\a{{\frak a}_{\mathbb V}}
\def\d{{\bf d}}
\def\e{{\bf e}}
\def\k{{\bf k}}
\def\n{{\bf n}}
\def\m{{\bf m}}
\def\p{{\bf p}}
\def\q{{\bf q}}
\def\x{{\bf x}}
\def\y{{\bf y}}
\def\To{{\longrightarrow}}
\def\aa{{\frak a}_B}
\def\sign{\operatorname{sgn}}
\def\1{{\bf 1}}
\def\0{{\bf 0}}
\def\r{{\mathbf r}^{\mathbb L}}
\def\mba{{\bf a}}
\def\dd{{\frak d}}


\title[Multigraded regularity and $a^*$-invariant]{Multigraded regularity, $a^*$-invariant and the minimal free resolution}
\author{Huy T\`ai H\`a}
\address{Tulane University, Department of Mathematics, 6823 St. Charles Ave., New Orleans LA 70118, USA}
\email{tai@math.tulane.edu}
\urladdr{http://www.math.tulane.edu/$\sim$tai/}

\begin{abstract}
In recent years, two different multigraded variants of Castelnuovo-Mumford regularity have been developed, namely {\it multigraded regularity}, defined by the vanishing of multigraded pieces of local cohomology modules, and the {\it resolution regularity vector}, defined by the multi-degrees in a minimal free resolution. In this paper, we study the relationship between multigraded regularity and the resolution regularity vector. Our method is to investigate multigraded variants of the usual $a^*$-invariant. This, in particular, provides an effective approach to examining the vanishing of multigraded pieces of local cohomology modules with respect to different multigraded ideals.
\end{abstract}
\maketitle



\section{Introduction}

\subsection{Objective.} Castelnuovo-Mumford regularity is a fundamental invariant in commutative algebra and algebraic geometry. Broadly speaking, regularity measures the computational complexity of a module or a sheaf. It is well known (cf. \cite{eg}) that  Castelnuovo-Mumford regularity can be defined either by the vanishing of graded pieces of local cohomology modules (equivalently, the vanishing of sheaf cohomology groups) or via bounds on the shifts in a graded minimal free resolution. These different ways of defining Castelnuovo-Mumford regularity are intrinsically linked and complementary to one another.

Many authors \cite{acn,hs,hv,hh,hh2,ms1,ms2,romer,sv,sv2} in recent years have developed multigraded variants of Castelnuovo-Mumford regularity. Inspired by the usual $\ZZ$-graded case, two different approaches have been investigated.

Maclagan and Smith \cite{ms1,ms2} (see also \cite{hh} for the bi-graded case) use local cohomology to define multigraded regularity. We recapture their definition for $\ZZ^k$-graded modules as follows (in fact, their definition works for $G$-graded modules where $G$ is a finitely generated abelian group). Let $S$ be a standard $\NN^k$-graded polynomial ring over a field $\k$, i.e. $S$ is the homogeneous coordinate ring of a multi-projective space $X = \PP^{N_1} \times \dots \times \PP^{N_k}$ for some $N_1, \dots, N_k \ge 1$. Let $\B$ be the $\NN^k$-graded irrelevant ideal of $S$, and let $F$ be a finitely generated $\ZZ^k$-graded $S$-module. For an integer $m$, let $\sign(m)$ be the sign function of $m$, and let $[m]_k$ be the finite set $\big\{ \p \in \NN^k \ \big| \ |\p| = |m| \big\}$, where $|\p|$ denotes the sum of coordinates of $\p$ and $|m|$ is the absolute value of $m$.

\begin{defn} \label{regularity-def}
Suppose $\m \in \ZZ^k$, then $F$ is said to be {\it $\m$-regular} if
$$H^i_{\B}(F)_{\n+\sign(1-i)\p} = 0 \ \text{for all} \ i \ge 0, \n \ge \m \ \text{and} \ \p \in \kset{1-i}_k.$$
The {\it multigraded regularity} of $F$ is defined to be the set
$$\Reg(F) = \big\{ \m \ \big| \ F \ \text{is} \ \m\text{-regular} \ \big\} \subseteq \ZZ^k.$$
\end{defn}

\begin{remark} \label{non-empty}
It is well known (cf. \cite[Theorem 1.6]{hy}) that $H^i_{\B}(F)_{\m} = 0$ for $\m \gg \0$. Thus, $\Reg(F) \not= \emptyset$ for any $\ZZ^k$-graded $S$-module $F$.
\end{remark}

The Serre-Grothendieck correspondence between sheaf cohomology on $X$ and local cohomology with respect to $\B$ gives corresponding notions of $\m$-regular and multigraded regularity for coherent $\O_X$-modules. Multigraded regularity shares similarities with the original definition of Castelnuovo-Mumford regularity (\cite{mum}), for instance, if a coherent sheaf $\F$ is $\m$-regular then $\F(\m)$ is generated by its global sections and the natural map $H^0(X, \F(\m)) \otimes H^0(X, \O_X(\n)) \rightarrow H^0(X, \F(\m+\n))$ is surjective for all $\n \in \NN^k$ (see \cite[Theorem 1.4]{ms1}).

On the other hand, developed from the theory of Hilbert functions and minimal free resolutions, Aramova, Crona and De Negri \cite{acn} (for $\NN^2$-graded), and Sidman and Van Tuyl \cite{sv} (for $\NN^k$-graded in general) extend the notion of Castelnuovo-Mumford regularity to the following notion of a resolution regularity vector.

\begin{defn} \label{resolution-reg-def}
For each $1 \le l \le k$, let $\e_l$ be the $l$th standard basis vector of $\ZZ^k$, and let
$$\rreg{l}(F) = \max \big\{ n_l \ \big| \ \exists i, \n = (n_1, \dots, n_k) \in \ZZ^k \ \text{so that} \ \Tor^S_i(F, \k)_{\n+i\e_l} \not= 0 \}.$$
The {\it resolution regularity vector} of $F$ is defined to be the vector
$$\Rreg(F) = (\rreg{1}(F), \dots, \rreg{k}(F)) \in \ZZ^k.$$
\end{defn}

\begin{remark} \label{resolution-def}
Let $\GG: 0 \rightarrow G_s \rightarrow \dots \rightarrow G_0 \rightarrow F \rightarrow 0$ be a minimal $\ZZ^k$-graded free resolution of $F$ over $S$, where $G_i = \bigoplus_j S(-\c{1}_{ij}, \dots, -\c{k}_{ij})$. We have $\Tor^S_i(F,\k) = H_i(\GG \otimes \k)$. Thus, if we set $\c{l}_i = \max_j \{ \c{l}_{ij} \}$, then it follows from the definition of $\rreg{l}(F)$ that
$$\rreg{l}(F) = \max_i \{ \c{l}_i - i \}.$$
We shall often make use of this fact.
\end{remark}

Unlike in the $\ZZ$-graded case where different techniques agree to give a unique invariant (see \cite{eg}), the relationship between multigraded regularity and the resolution regularity vector is not yet clear. The goal of this paper is to investigate this relationship.

The conceptual difficulty of the problem lies in the simple fact that bounded subsets of $\ZZ^k$ need not have single maximal or minimal elements. This makes it very hard to capture the vanishing or non-vanishing of multigraded pieces of local cohomology modules by a finite collection of vectors. Our work provides an effective approach for examining the vanishing of multigraded pieces of local cohomology modules with respect to different multigraded ideals. This allows us to pass from the $\ZZ^k$-graded case to the $\ZZ^r$-graded situation for $r < k$, or more generally, from finer gradings to coarser ones. Our results will not only give a bound for multigraded regularity in terms of the resolution regularity vector and other $\ZZ$-graded invariants but also exhibit an irregular behavior of the relationship between these invariants.

\subsection{Results and techniques.} We first observe that both multigraded regularity and the resolution regularity vector are characterized by the vanishing of multigraded pieces of local cohomology modules. In the $\ZZ$-graded case, the vanishing of graded pieces of local cohomology modules is controlled by the $a^*$-invariant. The relationship between Castelnuovo-Mumford regularity and $a^*$-invariant has also proved to be significant in many ways. Taking this as a starting point, we define two multigraded variants of the usual $\ZZ$-graded $a^*$-invariant, namely {\it multigraded $a^*$-invariant} $\aa^*(F) \subseteq \ZZ^k$ and the {\it $a^*$-invariant vector} $\a^*(F) \in \ZZ^k$. For $1 \le l \le k$, let $\M_l$ be the $S$-ideal generated by $S_{\e_l}$, the $l$th set of variables in $S$. Then multi-graded $a^*$-invariant and the $a^*$-invariant vector control the vanishing of multigraded pieces of local cohomology modules with respect to the $\NN^k$-graded irrelevant ideal $\B = \bigcap_{l=1}^k \M_l$ and the ideal $\M = \sum_{l=1}^k \M_l$, respectively. Our investigation is based on establishing a connection between multigraded $a^*$-invariant and the $a^*$-invariant vector and recovering a multigraded version of the well known relationship between regularity and $a^*$-invariant. Our techniques are inspired by previous work of the author and N.V. Trung \cite{ha1, ha2, ht}.

There are several spectral sequences relating local cohomology modules with respect to $\B$ and those with respect to $\M$. However, tracing through these spectral sequences to get any conclusion is a difficult task (especially when $k$ is large). Our method in establishing a connection between multigraded $a^*$-invariant and the $a^*$-invariant vector is to look at projection maps arising when $S$ is viewed with various multigraded structures, and examine how sheaf cohomology (and equivalently, local cohomology) behaves through a push forward along these morphisms. These projection maps can be thought of as the fiber spaces of the toric variety $X$ associated to faces of its nef cone. They give rise to a number of Mayer-Vietoris sequences relating local cohomology modules with respect to different multigraded ideals. Our first result, Theorem \ref{a-inv-lemma}, proves the vanishing of multigraded pieces of local cohomology modules with respect to different multigraded ideals when the multi-degrees are larger than the $a^*$-invariant vector. We show that for any $i \ge 0$, $\n \ge \a^*(F) + \1$ and any two disjoint index sets $I, J \subseteq \{ 1, \dots, k \}$, we have
$$H^i_{\M_{I,J}}(F)_\n = 0$$
where $\M_{I,J} = \big(\bigcap_{i \in I} \M_i\big) \bigcap \big( \sum_{j \in J} \M_j \big)$. Notice that $\M_{I,J}$ is the multigraded irrelevant ideal of the multi-projective space $\prod_{i \in I} \PP^{N_i} \times \PP^{\sum_{j \in J} (N_j+1)-1}$. This result is interesting on its own, but it also plays an essential role in establishing a lower bound for multigraded $a^*$-invariant. In \cite{ht}, the author and N.V. Trung introduced the notion of {\it projective $a^*$-invariant} which governed when sheaf cohomology could be passed through a blowing up morphism. We generalize this notion to that of {\it fiber $a^*$-invariant}. For each $l = 1, \dots, k$, let
$$\ppi{l}: X \rightarrow \PP^{N_1} \times \dots \times \widehat{\PP^{N_l}} \times \dots \times \PP^{N_k}$$
be the projection that drops the $l$th coordinate of $X$. The map $\ppi{l}$ is the fiber space associated to a facet of the nef cone of $X$. The fiber $a^*$-invariant $\la{l}^*(F)$ of $F$ with respect to $\ppi{l}$ determines the vanishing of higher direct images of $\F(\m)$ through $\ppi{l}$ (where $\F$ is the coherent sheaf associated to $F$ on $X$). It turns out that this important invariant not only governs when sheaf cohomology can be passed through projection maps but also provides an upper bound for multigraded $a^*$-invariant and multigraded regularity. We obtain in Theorem \ref{a-inv-bound} a bound for multigraded $a^*$-invariant in terms of the $a^*$-invariant vector and fiber $a^*$-invariant. More precisely, we prove that
\begin{align}
\a^*(F) + \NN^k \subseteq \aa^*(F) \subseteq \La^*(F) + \NN^k, \label{ainclusion}
\end{align}
where $\La^*(F) = \big(\la{1}^*(F), \dots, \la{k}^*(F)\big)$.

The study of multigraded regularity and the resolution regularity vector is more subtle than that of multigraded $a^*$-invariant and the $a^*$-invariant vector. We begin by proving another vanishing theorem (Theorem \ref{2lemma}) for multigraded pieces of local cohomology modules with respect to different multigraded ideals. Our next result, Theorem \ref{2regbound}, then gives a lower bound and an upper bound for $\Reg(F)$ in terms of $\Rreg(F)$ and other $\ZZ$-graded invariants. More precisely, we give a vector $\r(F)$ (see Theorem \ref{2regbound} for a precise statement) such that if $\delta = \operatorname{proj-dim} F$ then
\begin{align}
\bigcup_{\q \in \kset{\delta}_k} \big( \Rreg(F)+\delta\1 - \q + \NN^k \big) \subseteq \Reg(F) \subseteq \r(F)+\NN^k. \label{inclusionthm}
\end{align}
In \cite{sv}, Sidman and Van Tuyl obtained a similar (and slightly stronger) lower bound than that of (\ref{inclusionthm}) (with $\delta$ replaced by $\min \{ \dim X + 1, \operatorname{proj-dim} F \}$). Our method is completely different though. This result can be viewed as the first step toward obtaining a solution to the open problem \cite[Problem 4.10]{ms1}.
When $F$ is the coordinate ring of a set of fat points, Sidman and Van Tuyl improved their bound by proving that, in this situation, $\Rreg(F) + \NN^k \subseteq \Reg(F)$. They also raised the natural question of under what extra conditions is this inclusion still valid \cite[Question 2.6]{sv}. We give an example (Example \ref{resnotinreg}) showing that the inclusion $\Rreg(F) + \NN^k \subseteq \Reg(F)$ does not always hold. In our example, we in fact have the reverse inclusion $\Reg(F) \subsetneq \Rreg(F) + \NN^k$ (see Remark \ref{strictnoteq}). This shows that unlike multigraded $a^*$-invariant and the $a^*$-invariant vector (whose relationship is expressed as in (\ref{ainclusion})), multigraded regularity and the resolution regularity vector are not comparable in general.

We conclude the paper by showing that, just like multigraded regularity, the resolution regularity vector shares similarities with the original definition of Castelnuovo-Mumford regularity given in \cite{mum}. Beside the fact that $\Rreg(F)$ provides a crude bound for multi-degrees of syzygies in a minimal $\ZZ^k$-graded free resolution of $F$, we prove in Theorem \ref{globalgen} that if $\m \ge \Rreg(F)$ then $\F(\m)$ is generated by its global sections and the natural map $H^0(X, \F(\m)) \otimes H^0(X, \O_X(\n)) \rightarrow H^0(X, \F(\m+\n))$ is surjective for any $\n \ge \0$. This justifies the terminology resolution {\it regularity} vector.

Note that if $R$ is a standard $\NN^k$-graded $\k$-algebra and $F$ a $\ZZ^k$-graded $R$-module, then by considering the natural surjection $S \rightarrow R$ from a polynomial ring to $R$, we can view $F$ as a $\ZZ^k$-graded module over a standard $\NN^k$-graded polynomial ring. All our results, therefore, hold for any finitely generated $\ZZ^k$-graded module over a standard $\NN^k$-graded $\k$-algebra. All our techniques and results are independent of the ground field.

\subsection{Outlines of the paper.} In the next section, we give notation and terminology that will be used throughout the paper. In Section \ref{a-invariant-sec}, we introduce and study multigraded $a^*$-invariant, the $a^*$-invariant vector, fiber $a^*$-invariant, and the connection between these invariants. We start out by introducing multigraded $a^*$-invariant and recovering a multigraded version of the correspondence between regularity and $a^*$-invariant in $\S$\ref{ainvsubsec}. Next, in $\S$\ref{ainvvectsubsec}, we introduce the $a^*$-invariant vector and prove our first main theorem which establishes the vanishing of multigraded pieces of local cohomology modules with respect to different multigraded ideals. We continue in $\S$\ref{localainvsubsec} by introducing fiber $a^*$-invariant and conclude the section in $\S$\ref{ainvboundsubsec} by proving our next main theorem which relates multigraded $a^*$-invariant, the $a^*$-invariant vector and fiber $a^*$-invariant. In Section \ref{regularity-sec}, we investigate the relationship between multigraded regularity and resolution regularity. We first prove another vanishing theorem for multigraded pieces of local cohomology modules with respect to different multigraded ideals in $\S$\ref{vanishing-subsec}. The next part of the paper, $\S$\ref{reg-subsec}, is devoted to finding bounds on multigraded regularity in terms of the resolution regularity vector and other $\ZZ$-graded invariants. We conclude the paper in $\S$\ref{res-subsec} by showing that the resolution regularity vector has a Castelnuovo-type property on the global generation of twisted coherent sheaves.

\subsection{Acknowledgement.} The author thanks Dale Cutkosky for suggesting that Theorem \ref{globalgen} should be true, and thanks David Cox for observing that projection maps used in the paper correspond to fiber spaces of a toric variety. The author would also like to thank Adam Van Tuyl for fruitful discussions, and thank an anonymous referee for a careful reading and useful comments.


\section{Preliminaries} \label{notation-sec}

In this section we will introduce notation and terminology used in the paper. We refer the reader to \cite{bs, gd, har} for elementary facts about schemes, graded rings and local cohomology modules.

Throughout the paper, $\k$ will denote an arbitrary infinite field. Let $\ZZ$ be the set of integers, and let $\NN$ be the set of natural numbers including $0$. As usual, we will use $\e_1, \dots, \e_k$ to denote the standard basis vectors for $\ZZ^k$, and let $\0 = (0, \dots, 0)$ and $\1 = (1, \dots, 1)$. For $\m = (m_1, \dots, m_k)$ and $\n = (n_1, \dots, n_k)$ in $\ZZ^k$, we shall write $\m \ge \n$ if $m_i \ge n_i$ for all $i=1, \dots, k$. We write $\m \gg \0$ if $m_i \gg 0$ for all $i$. For a set $A \subseteq \ZZ^k$ and a vector $\n \in \ZZ^k$, let $A + \n = \big\{ \n+\m \ \big| \ \m \in A \big\}.$ Observe that $\n + \NN^k$ is the set $\big\{ \m \in \ZZ^k \ \big| \ \m \ge \n \big\}$. Recall that for an integer $m$,
$\kset{m}_k = \big\{ \n \in \NN^k \ \big| \ |\n| = |m| \big\}$
is a finite set.

Let $T$ be a standard $\NN^k$-graded $\k$-algebra (by {\it standard} we mean that $T$ is finitely generated over $T_\0 = \k$ by the elements of $T_{\e_1}, \dots, T_{\e_k}$) and let $\U \subseteq T$ be an $\NN^k$-graded homogeneous ideal. Then the local cohomology modules $H^i_\U(\cdot)$ can be defined in the category of $\ZZ^k$-graded $T$-modules as usual. If $M$ is a finitely generated $\ZZ^k$-graded $T$-module, then the $H^i_\U(M)$'s are also $\ZZ^k$-graded $T$-modules. If we omit the multigraded structure then these $\ZZ^k$-graded local cohomology modules coincide with the usual $\ZZ$-graded ones.

Our techniques in this paper are based upon the fact that local cohomology modules behave well under a change of grading. To be more precise, let $\psi: \ZZ^k \rightarrow \ZZ^r$ be a group homomorphism such that $\psi(\NN^k) \subseteq \NN^r$, and let $M^\psi = \bigoplus_{\m \in \ZZ^r} \big( \bigoplus_{\psi(\n) = \m} M_\n \big)$ for any $\ZZ^k$-graded $T$-module $M$. Then $T^\psi$ is a $\NN^r$-graded ring and $M^\psi$ is a finitely generated $\ZZ^r$-graded $T^\psi$-module. We can think of $M^\psi$ as $M$ viewed with a different grading. It can be seen that $\big(H^i_\U(\cdot)\big)^\psi$ and $H^i_{\U^\psi}(\cdot^\psi)$ are both $\delta$-functors and coincide when $i = 0$. Thus, $H^i_\U(M)^\psi = H^i_{\U^\psi}(M^\psi)$. For this reason, we shall use $H^i_\U(M)^\psi$ and $H^i_{\U^\psi}(M^\psi)$ interchangeably. Passing through a group homomorphism $\psi: \ZZ^k \rightarrow \ZZ^r$ can also be viewed as coarsening the grading. This approach is later taken in \cite{hh2,sv2}.

Throughout the paper, $\phi$ will denote the group homomorphism $\ZZ^k \rightarrow \ZZ$ given by $\n \mapsto |\n|$. For a nonempty index set $I \subseteq \{1, \dots, k\}$, let $\phi_I$ and $\phi_{\widehat{I}}$ denote the group homomorphisms $\ZZ^k \rightarrow \ZZ^{|I|}$ and $\ZZ^k \rightarrow \ZZ^{k-|I|}$ given by $\n \mapsto (n_i)_{i \in I}$ and $\n \mapsto (n_i)_{i \not\in I}$, respectively. Let $\n_I = \phi_I(\n)$ and $\n_{\widehat{I}} = \phi_{\widehat{I}}(\n)$. For simplicity of notation, when $I = \{l\}$, we use $\phi_l, \phi_{\widehat{l}}, \n_l$ and $\n_{\widehat{l}}$ to denote $\phi_I, \phi_{\widehat{I}}, \n_I$ and $\n_{\widehat{I}}$, respectively. Notice that $\n_l = n_l$ is the $l$th coordinate of $\n$, and so we shall use $\n_l$ and $n_l$ interchangeably.

For an $\NN$-graded ring $T$, we shall always use $T_+ = \bigoplus_{n > 0}T_n$ to denote its $\NN$-graded irrelevant ideal. If $M$ is a finitely generated $\ZZ$-graded $T$-module, then it is well known that $H^i_{T_+}(M)_n = 0$ for all $i \ge 0$ and $n \gg 0$. For $i \ge 0$, let
$$a^i(M) = \left\{ \begin{array}{lcl} -\infty & \text{if} & H^i_{T_+}(M) = 0 \\
\max \big\{ n \big| H^i_{T_+}(M)_n \not= 0 \big\} & \text{if} & H^i_{T_+}(M) \not= 0. \end{array} \right.$$
The {\it $a^*$-invariant} of $M$ is defined to be $a^*(M) = \max_i \{ a^i(M) \}$. The $a^*$-invariant was first introduced by Goto and Watanabe \cite{gw}, and then further studied by Sharp \cite{sharp} and Trung \cite{trung}. We shall use the terminology {\it $a$-invariant} to address $a^i(M)$ and $a^*(M)$ all together. Note that originally the term $a$-invariant often only referred to $a^{\dim M}(M)$. It is well known that
\begin{align}
\reg(M) = \max_i \{ a^i(M) + i \}, \label{reg-ainv}
\end{align}
where $\reg(M)$ denotes the usual $\ZZ$-graded Castelnuovo-Mumford regularity of $M$. Let
$\GG: 0 \rightarrow G_s \rightarrow \dots \rightarrow G_0 \rightarrow M \rightarrow 0$
be a minimal free resolution of $M$ over $T$, where $G_i = \bigoplus_j T(-c_{ij})$, and let $c_i = \max_j \{ c_{ij} \}$ for $i = 0, \dots, s$. Then $\reg(M)$ has another well known interpretation via $\GG$, namely
$\reg(M) = \max_i \{ c_i - i \}.$

From now on, let $S = \bigoplus_{\m \in \NN^k} S_\m$ be the standard $\NN^k$-graded coordinate ring of a multi-projective space $X = \PP^{N_1} \times \dots \times \PP^{N_k}$ over $\k$ (for some $k \ge 2$).
Recall that $\M_i$ denotes the $S$-ideal generated by $S_{\e_i}$, $\M = \sum_{i=1}^k \M_i$ and $\B = \bigcap_{i = 1}^k \M_i$. Recall further that for any two disjoint index sets $I,J \subseteq \{1, \dots, k\}$, $\M_I = \sum_{i \in I}\M_i$ and $\M_{I,J} = \big(\bigcap_{i \in I}\M_i\big) \bigcap \big(\sum_{j \in J}\M_j\big)$. Note that $\M_I = \M_{\emptyset,I}$ and $\B = \M_{\{1, \dots, k\},\emptyset} = \M_{\{1, \dots, k-1\}, \{k\}}$.

Let $F = \bigoplus_{\n \in \ZZ^k} F_\n$ be a finitely generated $\ZZ^k$-graded $S$-module. For each $1 \le l \le k$, clearly $F^{\phi_l}$ is a finitely generated $\ZZ$-graded $S^{\phi_l}$-module. For $q \in \ZZ$, we shall denote by $\lF{l}_q$ the degree $q$ piece of $F^{\phi_l}$.
Let $\F$ be the associated coherent sheaf of $F$ on $X$. The multigraded Serre-Grothendieck correspondence gives the exact sequence
$$0 \rightarrow H^0_{\B}(F) \rightarrow F \rightarrow \bigoplus_{\n \in \ZZ^k}H^0(X, \F(\n)) \rightarrow H^1_{\B}(F) \rightarrow 0$$
and isomorphisms $\bigoplus_{\n \in \ZZ^k}H^i(X, \F(\n)) \simeq H^{i+1}_{\B}(F)$ for $i \ge 1$. 


\section{Local cohomology modules and multigraded $a^*$-invariant} \label{a-invariant-sec}

In this section, we shall introduce multigraded $a^*$-invariant, the $a^*$-invariant vector and fiber $a^*$-invariant, and establish a connection between these invariants. We shall also prove a vanishing theorem for multigraded pieces of local cohomology modules with respect to different multigraded ideals.


\subsection{Multigraded $a^*$-invariant.} \label{ainvsubsec} We start by introducing a multigraded variant of the usual $a^*$-invariant using the vanishing of multigraded pieces of local cohomology modules with respect to the $\NN^k$-graded irrelevant ideal $\B$, and relate this invariant to multigraded regularity. For simplicity, let $d_l = N_l+1$ for $l = 1, \dots, k$, and let $D = \dim X = N_1 + \dots + N_k$. Let $d = \sum_{i = 1}^k d_i = \dim S$ and let $\d = (d_1, \dots, d_k)$. 

\begin{defn} \label{a-invariant-def}
For $i \ge 0$, let
$$\aa^i(F) = \big\{ \mba-\1 \in \ZZ^k \big| H^i_{\B}(F)_{\n} = 0 \ \forall \ \n \ge \mba \big\}.$$
The {\it multigraded $a^*$-invariant} of $F$ is defined to be the set
$\aa^*(F) = \bigcap_i \aa^i(F).$
\end{defn}
Again, it is known that $H^i_{\B}(F)_{\m} = 0$ for $\m \gg \0$ (cf. \cite[Theorem 1.6]{hy}). Thus, $\aa^i(F) \not= \emptyset$ for all $i \ge 0$. Clearly, by the Serre-Grothendieck correspondence between sheaf cohomology on $X$ and local cohomology with respect to $\B$ and Serre's vanishing theorem, we have $\aa^i(F) = \ZZ^k$ for $i > D+1$, so $\aa^*(F) = \bigcap_{i=0}^{D+1} \aa^i(F)$ is a finite intersection.

\begin{remark} \label{defremark}
When $k=1$ our definition of the $a$-invariants of $F$ agrees with the usual $a$-invariants in the following sense: $\aa^i(F) = \big\{ n \in \ZZ \ \big| \ n \ge a_i(F) \big\}$ and $\aa^*(F) = \big\{ n \in \ZZ \ \big| \ n \ge a^*(F) \big\}$.
\end{remark}

\begin{example} \label{example1}
In this example we consider the simplest case where $F = S$. For $i > D+1$ it is already known that $\aa^i(S) = \ZZ^k$. We shall find $\aa^i(S)$ for $0 \le i \le D+1$. For each nonempty index set $I \subseteq \{ 1, \dots, k \}$, denote $N_I = \sum_{i \in I} N_i$. Let $\N$ be the set $\big\{ N_I \ \big| \ I \subseteq \{1, \dots, k \} \big\}$. By K\"unneth's formula, we have
$$H^i\big(X, \O_X(\n)\big) = \sum_{j_1 + \dots + j_k = i} H^{j_1}\big(\PP^{N_1}, \O_{\PP^{N_1}}(n_1)\big) \otimes \dots \otimes H^{j_k}\big(\PP^{N_k}, \O_{\PP^{N_k}}(n_k)\big).$$
Moreover, $H^j\big(\PP^{N_l}, \O_{\PP^{N_l}}(n_l)\big)$ can be computed explicitly (cf. \cite[Theorem III.5.1]{har}).
Thus, for $i = 0$ we have $H^0\big(X, \O_X(\n)\big) = S_{\n}$. This, together with the Serre-Grothendieck correspondence, implies that $\aa^0(S) = \aa^1(S) = \ZZ^k$. For $i > 0$ and $i \not\in \N$, we have $H^i\big(X, \O_X(\n)\big) = 0$ for any $\n \in \ZZ^k$ (since in this case, for any decomposition $i = j_1 + \dots + j_k$, there always exists an $l$ such that $0 < j_l < N_l$). This and the Serre-Grothendieck correspondence imply that $\aa^{i+1}(S) = \ZZ^k$ for any $i \ge 1$ and $i \not\in \N$. It remains to consider $i \in \N$. Let $\A_i = \big\{ I \subseteq \{ 1, \dots, k \} \big| \ N_I = i \big\}$. We say that a set $G \subseteq \{1, \dots, k\}$ is a {\it generator} for $\A_i$ if $I \cap G \not= \emptyset$ for any $I \in \A_i$. We say that $G$ is a {\it minimal generator} of $\A_i$ if $G$ is a generator of $\A_i$ and any proper subset of $G$ is not. It follows from K\"unneth's formula that $H^i\big(X, \O_X(\n)\big) \not= 0$ if and only if there exists $I \in \A_i$ so that $n_l < -N_l \ \forall \ l \in I$ and $n_t \ge 0 \ \forall \ t \not\in I$. Hence, by the Serre-Grothendieck correspondence we have $\aa^{i+1}(S) = \big\{ \n \in \ZZ^k \ \big| \ \exists \ \text{a minimal generator $G$ of} \ \A_i \ \text{such that} \ n_l \ge -N_l \ \forall \ l \in G \big\}$.

By considering $\aa^{N_l+1}(S)$ for $l = 1, \dots, k$, we also get $\aa^*(S) = \big\{ \n -\1 \in \ZZ^k \ \big| \ n_l \ge -N_l \ \forall \ l = 1, \dots, k \big\} = -\d+\NN^k$.
\end{example}

The following proposition gives a multigraded version of the relationship between regularity and $a$-invariant.

\begin{proposition} \label{2a-invariant}
Let $F$ be a finitely generated $\ZZ^k$-graded $S$-module. Then
$$\Reg(F) = \bigcap_{i=0}^{D+1} \bigcap_{\p \in \kset{1-i}_k} \big(\aa^i(F) +\1 -\sign(1-i)\p \big).$$
\end{proposition}

\begin{proof} By definition, $\m \in \Reg(F)$ if and only if
$H^i_{\B}(F)_{\n+\sign(1-i)\p} = 0 \ \forall \ i \ge 0, \n \ge \m, \p \in \kset{1-i}_k.$
In other words, $\m \in \Reg(F)$ if and only if $\m + \sign(1-i)\p -\1 \in \aa^i(F)$ for all $i \ge 0$ and $\p \in \kset{1-i}_k$. Thus, $\m \in \Reg(F)$ if and only if $\m \in \aa^i(F) +\1 - \sign(1-i)\p$ for all $i \ge 0$ and $\p \in \kset{1-i}_k$. Hence, $\m \in \Reg(F)$ if and only if
$$\m \in \bigcap_{i \ge 0} \bigcap_{\p \in \kset{i-1}_k} \big(\aa^i(F) +\1 - \sign(1-i)\p \big).$$
Since $\aa^i(F) = \ZZ^k$ for all $i > D+1$, the above intersection is, in fact, a finite intersection as $i$ runs from $0$ to $D+1$. The proposition is proved.
\end{proof}

\begin{example} \label{example2}
It follows from the calculation in Example \ref{example1}, by considering $\aa^{N_l+1}(S)$ for $l = 1, \dots, k$ and applying Proposition \ref{2a-invariant}, that $\Reg(S) = \NN^k$. From this we also have $\Reg\big(S(-\n)\big) = \n + \NN^k$.
\end{example}

\begin{remark} \label{regvsainv}
Proposition \ref{2a-invariant} implies that $\Reg(F) \subseteq \bigcap_{\p \in \kset{1}_k} \big( \aa^0(F) +\1 - \p\big) = \aa^0(F)+\1.$ Clearly, Proposition \ref{2a-invariant} also implies that $\Reg(F) \subseteq \bigcap_{\p \in \kset{1-i}_k} \big( \aa^i(F)+\1+\p \big) \subseteq \aa^i(F)+\1 $ for any $i > 0$. Hence, $\Reg(F) \subseteq \aa^*(F)+\1$.
\end{remark}


\subsection{The $a^*$-invariant vector and the vanishing of multigraded pieces of local cohomology modules.} \label{ainvvectsubsec} In $\S$\ref{ainvsubsec} we have defined the multigraded $a^*$-invariant which controls the vanishing of multigraded pieces of local cohomology modules with respect to $\B$. In this subsection we shall define the $a^*$-invariant vector which governs the vanishing of multigraded pieces of local cohomology modules with respect to $\M$. We shall also prove our first main result, a vanishing theorem for local cohomology modules with respect to different multigraded ideals.

We start with a lemma which allows us to pass between local cohomology modules with respect to the $\NN$-graded irrelevant ideal and the maximal homogeneous ideal.

\begin{lemma} \label{hyrylemma}
Let $S$ be a standard $\NN$-graded algebra over $R$, where $R$ is either a local ring or a standard $\NN$-graded $\k$-algebra. Let $\mm$ be the maximal ideal of $R$ if $R$ is a local ring, or the maximal homogeneous ideal of $R$ if $R$ is an $\NN$-graded $\k$-algebra. Let $S_+$ be the $\NN$-graded irrelevant ideal of $S$ and let $\M = \mm \bigoplus S_+$. Let $F$ be a finitely generated $\ZZ$-graded $S$-module and $n_0 \in \ZZ$. Then $H^i_\M(F)_n = 0$ for all $i \ge 0$ and $n \ge n_0$ if and only if $H^i_{S_+}(F)_n = 0$ for all $i \ge 0$ and $n \ge n_0$.
\end{lemma}

\begin{proof} By moving to the completion $\widehat{S} = S \otimes_R \widehat{R}$, where $\widehat{R}$ is the completion of $R$ with respect to $\mm$, we may assume that $R$ is a complete local ring with maximal ideal $\mm$. The conclusion now follows from that of \cite[Lemma 2.3]{hy}.
\end{proof}

\begin{defn} \label{a-invariant-def2}
For each $l = 1, \dots, k$, let
$$a_{l,i}(F) = \left\{ \begin{array}{lcl} -\infty & \text{if} & H^i_\M(F) = 0 \\
\max \big\{ n_l \ \big| \ \exists \n = (n_1, \dots, n_k) \in \ZZ^k : H^i_\M(F)_\n \not= 0 \big\} & \text{if} & H^i_\M(F) \not= 0 \end{array} \right.$$
for $i \ge 0$, and let
$a_{l,*}(F) = \max_i \big\{ a_{l,i}(F) \big\}.$
The {\it $a^*$-invariant vector} of $F$ is defined to be the vector
$$\a^*(F) = \big(a_{1,*}(F), \dots, a_{k,*}(F)\big) \in \ZZ^k.$$
\end{defn}

\begin{remark} Observe that $S^{\phi_l}$ is a standard $\NN$-graded ring with the $\NN$-graded irrelevant ideal $\M_l$ and $F^{\phi_l}$ is a finitely generated $\ZZ$-graded $S^{\phi_l}$-module. Thus, $H^i_{\M_l}(F)^{\phi_l}$ is a $\ZZ$-graded $S^{\phi_l}$-module and $[H^i_{\M_l}(F)^{\phi_l}]_q = 0$ for $q \gg 0$. Moreover, $H^i_{\M_l}(F)^{\phi_l} = 0$ for $i \gg 0$. Therefore, by Lemma \ref{hyrylemma}, $[H^i_\M(F)^{\phi_l}]_q = 0$ for all $i \ge 0$ and $q \gg 0$. Hence, in Definition \ref{a-invariant-def2}, $a_{l,i}(F) < \infty$ for all $i$ and $a_{l,*}(F)$ is well-defined.
When $k = 1$, the $a^*$-invariant vector agrees with the usual $a^*$-invariant, namely $\a^*(F) = a^*(F) \in \ZZ$.
\end{remark}

\begin{example} \label{example3}
Consider $F = S$. As mentioned in Section \ref{notation-sec}, local cohomology functors can be defined on the category of $\ZZ^k$-graded $S$-modules. Thus, by considering the \v{C}ech complex associated to the elements of $\bigcup_{i=1}^k S_{\e_i}$ and keeping track of the multi-degrees, we have a perfect pairing $S_\n \times H^d_\M(S)_{-\n - \d} \rightarrow \k$, and $H^i_\M(S) = 0$ for any $i < d$. This implies that $\a^*(S) = -\d$.
\end{example}






The next lemma is crucial in our proof of the first main result.

\begin{lemma} \label{a-inv-lem1}
Let $S$ be a standard $\NN^k$-graded polynomial ring over a field $\k$ and let $F$ be a finitely generated $\ZZ^k$-graded $S$-module. Then, for any index set $I \subseteq \{ 1, \dots, k\}$, any $i \ge 0$ and $\n \in \ZZ^k$ such that $|\n_I| > \sum_{l \in I} a_{l,*}(F)$, we have
$H^i_{\M_I}(F)_\n = 0.$
\end{lemma}

\begin{proof} If $I = \emptyset$ then the statement is vacuously true. Assume that $I \not= \emptyset$. Let $\varphi: \ZZ^{|I|} \rightarrow \ZZ$ be the group homomorphism given by $\m \mapsto |\m|$. Let $\psi = \varphi \circ \phi_I$. Observe that $S^\psi$ is a standard $\NN$-graded ring with the $\NN$-graded irrelevant ideal $\M_I$. Moreover, $H^i_\M(F)^\psi$ and $H^i_{\M_I}(F)^\psi$ are $\ZZ$-graded $S^\psi$-modules.

Let $m = \sum_{l \in I} a_{l,*}(F)$. It can be seen that for any $\n \in \ZZ^k$ with $|\n_I| > m$ there exists $j \in I$ such that $n_j > a_{j,*}(F)$, which then implies that $H^i_\M(F)_\n = 0$ for any $i \ge 0$. Thus, $[H^i_\M(F)^\psi]_n = 0$ for all $i \ge 0$ and $n > m$. This and Lemma \ref{hyrylemma} imply that $[H^i_{\M_I}(F)^\psi]_n = 0$ for any $i \ge 0$ and $n > m$. That is, for any $\n \in \ZZ^k$ with $|\n_I| > \sum_{l \in I} a_{l,*}(F)$ and $i \ge 0$, we have $H^i_{\M_I}(F)_\n = 0$. The lemma is proved.
\end{proof}

We are now ready to prove a vanishing theorem for multigraded pieces of local cohomology modules with respect to different multigraded ideals.

\begin{theorem} \label{a-inv-lemma}
Let $S$ be a standard $\NN^k$-graded polynomial ring over a field $\k$ and let $F$ be a finitely generated $\ZZ^k$-graded $S$-module. Then, for any $\n \ge \a^*(F) + \1, i \ge 0$ and any two disjoint index sets $I, J \subseteq \{ 1, \dots, k \}$, we have
\begin{align}
H^i_{\M_{I,J}}(F)_\n = 0. \label{a-inv-lemma-eq}
\end{align}
\end{theorem}

\begin{proof} We shall use induction on $r = |I \cup J|$, the total number of elements in $I$ and $J$. If $r = 0$ then the statement is vacuous. If $r = 1$ then the statement clearly follows from Lemma \ref{a-inv-lem1}. Suppose $r \ge 2$. Without loss of generality, assume that $I = \{ 1, \dots, s \}$ and $J = \{s+1, \dots, r\}$. If $s = 0$ then the statement again follows from Lemma \ref{a-inv-lem1}. Suppose $s > 0$. By identifying $\M_{I, \emptyset}$ with $\M_{I \backslash \{ s \}, \{s\}}$ we may assume that $J \not= \emptyset$, i.e. $s < r$. We now use ascending induction on $s$. Let $K = I \backslash \{s\}$ and $L = J \cup \{s\}$. Observe that $\M_{K,J} + \M_{I, \emptyset} = \M_{K,L}$ and $\M_{K,J} \cap \M_{I, \emptyset} = \M_{I,J}$. Consider the following Mayer-Vietoris sequence of local cohomology modules
$$\dots \rightarrow H^i_{\M_{K,J}}(F) \oplus H^i_{\M_{I, \emptyset}}(F) \rightarrow H^i_{\M_{I,J}}(F) \rightarrow H^{i+1}_{\M_{K,L}}(F) \rightarrow \dots $$
By induction on $s$, we have
$H^{i+1}_{\M_{K,L}}(F)_\n = 0$
for all $i \ge 0$ and $\n \ge \a^*(F)+\1$. By induction on $r$, we have
$H^i_{\M_{K,J}}(F)_\n = 0 = H^i_{\M_{I, \emptyset}}(F)_\n$
for all $i \ge 0$ and $\n \ge \a^*(F)+\1$. Hence,
$H^i_{\M_{I,J}}(F)_\n = 0$
for all $i \ge 0$ and $\n \ge \a^*(F)+\1$. The result is proved.
\end{proof}


\subsection{Fiber $a^*$-invariant.} \label{localainvsubsec} Our notion of fiber $a^*$-invariant is a generalization of that of {\it projective $a^*$-invariant} introduced by the author and N.V. Trung \cite{ht} in the study of the Cohen-Macaulay property of projective embeddings of blowup schemes.

Fix an integer $l \in \{ 1, \dots, k \}$. Let $\ppi{l}$ be the projection
$$\ppi{l}: X = \PP^{N_1} \times \dots \times \PP^{N_k} \To \PP^{N_l} \times \dots \times \widehat{\PP^{N_l}} \times \dots \times \PP^{N_k}$$
and let $\X{l} = \ppi{l}(X)$. Recall that $S^{\phi_l}$ is a standard $\NN$-graded ring and $F^{\phi_l}$ is a finitely generated $\ZZ$-graded $S^{\phi_l}$-module. Let $S^{\phi_l}_+$ denote the $\NN$-graded irrelevant ideal of $S^{\phi_l}$. Let $\mS{l} = [S^{\phi_l}]_0$ be the degree $0$ piece of $S^{\phi_l}$. It is easy to see that $\mS{l}$ is the $\NN^{k-1}$-graded homogeneous coordinate ring of $\X{l}$ and, for $q \in \ZZ$, $\lF{l}_q = [F^{\phi_l}]_q$ is a $\ZZ^{k-1}$-graded $\mS{l}$-module. Let $\widetilde{\lF{l}_q}$ denote the coherent sheaf associated to $\lF{l}_q$ on $\X{l}$.

Recall further that $S^{\phi_{\widehat{l}}}$ is a standard $\NN^{k-1}$-graded ring and $F^{\phi_{\widehat{l}}}$ is a finitely generated $\ZZ^{k-1}$-graded $S^{\phi_{\widehat{l}}}$-module.
By composing with the natural inclusion $\mS{l} \hookrightarrow S$ and the isomorphism $S \simeq S^{\phi_{\widehat{l}}}$, we can view $F^{\phi_{\widehat{l}}}$ as a $\ZZ^{k-1}$-graded $\mS{l}$-module.
For each $\ZZ^{k-1}$-graded homogeneous prime ideal $\pp \subseteq \mS{l}$, let $\mS{l}_{\pp}$ be the localization of $\mS{l}$ at $\pp$ and let $F_\pp = F^{\phi_{\widehat{l}}} \otimes_{\mS{l}} \mS{l}_\pp$. Then, $F_\pp$ is a $\ZZ^{k-1}$-graded $S_\pp$-module. Denote by $F_{(\pp)}$ the collection of all elements of degree $\0 \in \ZZ^{k-1}$ of $F_\pp$, and call this the {\it homogeneous localization} of $F$ at $\pp$. In this situation, homogeneous localization at any $\NN^{k-1}$-graded prime $\pp \subseteq \mS{l}$ annihilates the $\ZZ^{k-1}$-graded structure of an $S^{\phi_{\widehat{l}}}$-module. Incorporating this with the original $\ZZ^k$-graded structure of $F$, it can be seen that $F_{(\pp)}$ is a $\ZZ$-graded $S_{(\pp)}$-module. Thus, the usual $\ZZ$-graded $a$-invariants are defined, namely
$a^i(F_{(\pp)}) = \max \big\{ m \big| [H^i_{S_{(\pp)+}}(F_{(\pp)})]_m \not= 0 \big\}$ for $i \ge 0$.

\begin{defn} \label{local-a-invariant-def}
For each $l = 1, \dots, k$, let
$$\la{l}^i(F) = \max \big\{ a^i(F_{(\pp)}) \big| \pp \text{ an $\NN^{k-1}$-graded prime ideal in } \mS{l} \big\} \ \text{for} \ i \ge 0.$$
The $l$th {\it fiber $a^*$-invariant} of $F$ is defined to be
$$\la{l}^*(F) = \max_i \big\{ \la{l}^i(F) \big\}.$$
\end{defn}

Note that $H^i_{S_{(\pp)+}}(F_{(\pp)}) = \big(H^i_{S^{\phi_{l}}_+}(F^{\phi_l})\big)_{(\pp)}$ (cf. \cite[Remark 2.2]{sharp}). Thus, $\la{l}^i(F) \le a^i(F^{\phi_l})$ is a finite number for any $i \ge 0$. Note also that since $\dim S_{(\pp)} \le \dim S = d$, we have $\la{l}^i(F) = 0$ for all $i > d$. Thus, $\la{l}^*(F)$ is well defined as the maximum of finitely many numbers.

For $\m = (m_1, \dots, m_k) \in \ZZ^k$ and each $l \in \{ 1, \dots, k \}$, we shall denote by $\m^{[l]}$ the vector $(m_1, \dots, m_{l-1}, \la{l}^*(F), m_{l+1}, \dots, m_k) \in \ZZ^k$. The fiber $a^*$-invariant of $F$ provides measures of when one can pass from sheaf cohomology on $X$ to that on its projection $\X{l}$. The following statement is a generalization of \cite[Proposition 1.3 and Theorem 2.3]{ht}.

\begin{proposition} \label{2projection}
Let $\F$ be the coherent sheaf on $X$ associated to the $\ZZ^k$-graded $S$-module $F$. Then for each $l = 1, \dots, k$, we have
\begin{enumerate}
\item $\ppi{l}_* \F(\m) = \widetilde{\lF{l}_{m_l}}(\m_{\widehat{l}}) \ \text{and} \ R^j \ppi{l}_* \F(\m) = 0$
for any $j > 0$ and $\m \in \ZZ^k$ such that $m_l > \la{l}^*(F)$.
\item $\m^{[l]} -\1 \not\in \aa^*(F)$ for any $\m \in \ZZ^k$ such that $\m_{\widehat{l}} \gg \0$.
\end{enumerate}
\end{proposition}

\begin{proof} We shall write $\F^{\phi_l}(q)$ to denote the twisted $\O_X$-modules with respect to this $\ZZ$-graded structure of $F^{\phi_l}$. Then,
$\F(\m) = \F^{\phi_l}(m_l) \otimes {\ppi{l}}^* \O_{\X{l}}(\m_{\widehat{l}}).$
Thus, by the projection formula,
\begin{align}
\ppi{l}_* \F(\m) = \ppi{l}_* \F^{\phi_l}(m_l) \otimes \O_{\X{l}}(\m_{\widehat{l}}) \text{ and }
R^j \ppi{l}_* \F(\m) = R^j \ppi{l}_* \F^{\phi_l}(m_l) \otimes \O_{\X{l}}(\m_{\widehat{l}}). \label{2projection-eq100}
\end{align}

Let $\pp$ be any closed point of $\X{l}$, and consider the restriction $\ppi{l}_\pp$ of $\ppi{l}$ over an affine open neighborhood $\spec \O_{\X{l}, \pp}$ of $\pp$,
$\ppi{l}_\pp: X_\pp = X \times_{\X{l}} \spec \O_{\X{l}, \pp} \To \spec \O_{\X{l}, \pp}.$
We have $\F\big|_{X_\pp} = \widetilde{F_{(\pp)}}$ where $\widetilde{F_{(\pp)}}$ is the sheaf associated to $F_{(\pp)}$ on $X_\pp$. Therefore, 
\begin{align}
R^j \ppi{l}_* \F^{\phi_l}(m_l) \big|_{\spec \O_{\X{l},\pp}} & = R^j {\ppi{l}_\pp}_* \big(\widetilde{F_{(\pp)}}(m_l)\big) = H^j(X_\pp, \widetilde{F_{(\pp)}}(m_l))\widetilde{\ } \ \ \forall \ j \ge 0. \label{2proj0}
\end{align}
By the Serre-Grothendieck correspondence we have an exact sequence
\begin{align}
0 \To H^0_{S_{(\pp)+}}(F_{(\pp)})_{m_l} \To \big[F_{(\pp)}\big]_{m_l} \To H^0(X_\pp, \widetilde{F_{(\pp)}}(m_l)) \To H^1_{S_{(\pp)+}}(F_{(\pp)})_{m_l} \To 0 \label{2exact}
\end{align}
and isomorphisms
\begin{align}
H^i(X_\pp, \widetilde{F_{(\pp)}}(m_l)) \simeq H^{i+1}_{S_{(\pp)+}}(F_{(\pp)})_{m_l} \text{ for } i > 0. \label{2iso}
\end{align}

Suppose that $m_l > \la{l}^*(F)$. Then for any $j \ge 0$ and $\pp \in \X{l}$, we have $m_l > a_j(F_{(\pp)})$ and
$H^j_{S_{(\pp)+}}(F_{(\pp)})_{m_l} = 0.$
Thus, for any $j \ge 0$ and $\pp \in \X{l}$, (\ref{2proj0}) gives
$$R^j \ppi{l}_* \F^{\phi_l}(m_l) \big|_{\spec \O_{\X{l},\pp}} = H^j(X_\pp, \widetilde{F_{(\pp)}}(m_l))\widetilde{\ } \ = \left\{ \begin{array}{lcr} \widetilde{\big(\lF{l}_{m_l}\big)_{(\pp)}} & \text{for} & j = 0 \\ 0 & \text{for} & j > 0. \end{array} \right.$$
Hence, $\ppi{l}_* \F^{\phi_l}(m_l) = \widetilde{\lF{l}_{m_l}} \ \text{and} \ R^j \ppi{l}_* \F^{\phi_l}(m_l) = 0 \ \text{for} \ j > 0$. (1) now follows from (\ref{2projection-eq100}). 

For simplicity, we write $\la{l}$ for $\la{l}^*(F)$. Let $\lr{l}$ be the least integer such that $\la{l}^{\lr{l}}(F) = \la{l}^*(F)$. That is,
\begin{align}
\left\{ \begin{array}{rcl} H^i_{S_{(\pp)+}}(F_{(\pp)})_{\la{l}} & = 0 & \text{ for } i < \lr{l} \text{ and any } \pp \in \X{l}, \\
H^{\lr{l}}_{S_{(\qq)+}}(F_{(\qq)})_{\la{l}} & \not= 0 & \text{ for some } \qq \in \X{l}. \end{array} \right. \label{non-vanish}
\end{align}

If $\lr{l} \le 1$ then (\ref{2exact}) and (\ref{non-vanish}) give us $H^0(X_\qq, \widetilde{F_{(\qq)}}(\la{l})) \not= \big[F_{(\qq)}\big]_{\la{l}}$. This and (\ref{2proj0}) imply that
\begin{align}
\ppi{l}_* \F(\m^{[l]}) \not= \widetilde{\lF{l}_{\la{l}}}(\m_{\widehat{l}}) \text{ for all } \m_{\widehat{l}} \in \ZZ^{k-1}. \label{2noneq}
\end{align}
By Remark \ref{non-empty}, we have $\Reg \G \not= \emptyset$ for any coherent sheaf $\G$ on $\X{l}$. Thus, it follows from (\ref{2projection-eq100}) and \cite[Theorem 1.4]{ms1} that $\ppi{l}_* \F(\m^{[l]})$ and $\widetilde{\lF{l}_{\la{l}}}(\m_{\widehat{l}})$ are generated by global sections for $\m_{\widehat{l}} \gg \0$. Therefore, by (\ref{2noneq}) we have
$H^0(\X{l}, \ppi{l}_* \F(\m^{[l]})) \not= H^0(\X{l}, \widetilde{\lF{l}_{\la{l}}}(\m_{\widehat{l}})) = F_{\m^{[l]}}$
for $\m_{\widehat{l}} \gg \0$. Moreover
$H^0(X, \F(\m^{[l]})) = H^0(\X{l}, \ppi{l}_* \F(\m^{[l]})).$
Hence,
$$H^0(X, \F(\m^{[l]})) \not= F_{\m^{[l]}} \ \text{for} \ \m_{\widehat{l}} \gg \0.$$
This and the Serre-Grothendieck correspondence imply that for $\m_{\widehat{l}} \gg \0$, we have either $H^0_{\B}(F)_{\m^{[l]}} \not= 0$ or $H^1_{\B}(F)_{\m^{[l]}} \not= 0$. Thus, for $\m_{\widehat{l}} \gg \0$, we have either $\m^{[l]}-\1 \not\in \aa^0(F)$ or $\m^{[l]}-\1 \not\in \aa^1(F)$.

If $\lr{l} \ge 2$ then (\ref{non-vanish}), (\ref{2exact}) and (\ref{2iso}) imply that $H^i(X_\pp, \widetilde{F_{(\pp)}}(\la{l})) = 0$ for all $\pp \in \X{l}$, $0 < i < \lr{l}-1$, and $H^{\lr{l}-1}(X_\qq, \widetilde{F_{(\qq)}}(\la{l})) \not= 0$ for some $\qq \in \X{l}$. This and (\ref{2proj0}) imply that
\begin{align}
\left\{ \begin{array}{l} R^i \ppi{l}_* \F(\m^{[l]}) = 0 \ \text{for} \ 0 < i < \lr{l}-1, \\
R^{\lr{l}-1} \ppi{l}_* \F(\m^{[l]}) \not= 0.\end{array} \right. \label{high-img}
\end{align}
Moreover, by (\ref{2projection-eq100}) we have
$H^{\lr{l}-1}(\X{l}, \ppi{l}_* \F(\m^{[l]})) = 0 \ \text{for} \ \m_{\widehat{l}} \gg \0.$
Thus, using the Leray spectral sequence
$H^i(\X{l}, R^j \ppi{l}_* \F(\m^{[l]})) \Rightarrow H^{i+j}(X, \F(\m^{[l]}))$
and (\ref{high-img}), we can deduce that
$$H^{\lr{l}-1}(X, \F(\m^{[l]})) = H^0(\X{l}, R^{\lr{l}-1} \ppi{l}_* \F(\m^{[l]})) \text{ for } \m_{\widehat{l}} \gg \0.$$
Furthermore, $\Reg\big( R^{\lr{l}-1} \ppi{l}_* \F^{\phi_l}(\la{l})\big) \not= \emptyset$, so we have that $R^{\lr{l}-1} \ppi{l}_* \F(\m^{[l]})$ is generated by its global sections for $\m_{\widehat{l}} \gg \0$. Therefore, (\ref{high-img}) now gives
$$H^{\lr{l}-1}(X, \F(\m^{[l]})) \not= 0 \ \text{for} \ \m_{\widehat{l}} \gg \0.$$
This and the Serre-Grothendieck correspondence imply that $\m^{[l]} -\1 \not\in \aa^{\lr{l}}(F)$ for all $\m_{\widehat{l}} \gg \0$.

We have shown that, in any case, there exists $i \ge 0$ such that $\m^{[l]} -\1 \not\in \aa^i(F)$ for all $\m_{\widehat{l}} \gg \0$. The conclusion follows since $\aa^*(F) = \bigcap_{i \ge 0} \aa^i(F)$.
\end{proof}

\begin{corollary} \label{2projcor}
Suppose $\m \in \ZZ^k$ and $1 \le l \le k$. With the notation as in Proposition \ref{2projection}, we have
$$H^i(X, \F(\m)) = H^i(\X{l}, \ppi{l}_* \F(\m)) = H^i(\X{l}, \lF{l}_{m_l}(\m_{\widehat{l}}))$$
for any $i \ge 0$, $\m_{\widehat{l}} \in \ZZ^{k-1}$ and $m_l > \la{l}^*(F)$.
\end{corollary}

\begin{proof} The conclusion follows from part (1) of Proposition \ref{2projection} and the Leray spectral sequence $E^{i,j}_2 = H^i\big(\X{l}, R^j \ppi{l}_* \F(\m)\big) \Rightarrow H^{i+j}(X, \F(\m)).$
\end{proof}

\subsection{Bounding the multigraded $a^*$-invariant.} \label{ainvboundsubsec} We now establish a connection between multigraded $a^*$-invariant, the $a^*$-invariant vector and fiber $a^*$-invariant.

\begin{theorem} \label{a-inv-bound}
Let $S$ be a standard $\NN^k$-graded polynomial ring over a field $\k$ and let $F$ be a finitely generated $\ZZ^k$-graded $S$-module. Let $\La^*(F) = (\la{1}^*(F), \dots, \la{k}^*(F))$. Then,
\begin{align}
\a^*(F) + \NN^k \subseteq \aa^*(F) \subseteq \La^*(F) + \NN^k. \label{a-inv-bound-eq}
\end{align}
\end{theorem}

\begin{proof} Suppose $\n \ge \a^*(F)$. Pick $I = \{ 1, \dots, k \}$ and $J = \emptyset$. Then, $\M_{I,J} = \B$. By Theorem \ref{a-inv-lemma}, we have $H^i_{\B}(F)_\m = 0$ for any $i \ge 0$ and $\m \ge \n+\1$. This, by definition, implies that $\n \in \aa^*(F)$. The left hand inclusion of (\ref{a-inv-bound-eq}) is proved.

Suppose $\n \in \aa^*(F)$. Fix $l \in \{ 1, \dots, k \}$. Choose $\m \in \ZZ^{k-1}$ such that $m_l = n_l$ and $\m_{\widehat{l}} \gg \0 \in \ZZ^{k-1}$ (in particular, we shall choose $m_i \ge n_i+1$ for all $i \not= l$). If $n_l \le \la{l}^*(F) - 1$, then $\m^{[l]} - \1 \ge \n$, so $\m -\1 \in \aa^*(F)$ by definition. However, it follows from Proposition \ref{2projection} that $\m^{[l]} -\1 \not\in \aa^*(F)$ for all $\m_{\widehat{l}} \gg \0$. This contradiction implies that $n_l \ge \la{l}^*(F)$. Since this is true for any $l \in \{ 1, \dots, k \}$, we have $\n \in \La^*(F)+\NN^k$. The right hand inclusion of (\ref{a-inv-bound-eq}) is proved.
\end{proof}

\begin{example} \label{example4}
Consider $F = S$. As calculated in Example \ref{example3}, $\a^*(S) = -\d$. Moreover, it follows from Example \ref{example1} that $\aa^*(S) = -\d + \NN^k$. Hence, the left hand inclusion of (\ref{a-inv-bound-eq}) is an equality in this case. This shows that the left hand inclusion of (\ref{a-inv-bound-eq}) is sharp.

Observe further that for each $l \in \{ 1, \dots, k \}$ and $\pp \in \X{l}$, we have $S_{(\pp)} = S \otimes_{\mS{l}} \mS{l}_{(\pp)} = \mS{l}_{(\pp)}[S_{\e_l}]$ is a standard $\NN$-graded polynomial ring over $\mS{l}_{(\pp)}$. Thus, it is well known that $a^*\big(S_{(\pp)}\big) = -N_l-1 = -d_l$. This implies that $\la{l}^*(S) \ge -d_l$ for any $l = 1, \dots, k$. Therefore, $\La^*(S) \ge -\d$. Moreover, it follows from Theorem \ref{a-inv-bound} and Example \ref{example3} that $\La^*(S) \le \a^*(S) = -\d$. Hence, $\La^*(S) = -\d$ and the right hand inclusion of (\ref{a-inv-bound-eq}) is an equality in this case. This shows that the right hand inclusion of (\ref{a-inv-bound-eq}) is also sharp.
\end{example}


\section{Multigraded regularity and the minimal free resolution} \label{regularity-sec}

In this section we investigate multigraded regularity and the resolution regularity vector, and the relationship between these invariants and other $\ZZ$-graded invariants. 


\subsection{A vanishing theorem.} \label{vanishing-subsec} We start by proving a vanishing theorem for multigraded pieces of local cohomology modules with respect to different multigraded ideals. Our Theorem \ref{2lemma} serves as the analog of Theorem \ref{a-inv-lemma} in bounding the multigraded regularity. However, the proof of Theorem \ref{2lemma} is more subtle than that of Theorem \ref{a-inv-lemma} and requires more ideas. The technical difficulty lies handling in the shifts $\sign(1-i)\p$ in the definition of multigraded regularity.

\begin{lemma} \label{2lem00}
Let $S$ be a standard $\NN^2$-graded polynomial ring over a field $\k$ and let $F$ be a finitely generated $\ZZ^2$-graded $S$-module. Let $\mm_1 = \M_1 \cap \k[S_{(1,0)}]$. Then,
$$H^i_{\M_1}(F)_{(m,n)} = H^i_{\mm_1}\big(\lF{2}_n\big)_m.$$
\end{lemma}

\begin{proof} The proof follows from a direct application of that of \cite[Lemma 2.1]{cht}.
\end{proof}

\begin{lemma} \label{2lem0}
Let $S$ be a standard $\NN^k$-graded polynomial ring over a field $\k$ and let $F$ be a finitely generated $\ZZ^k$-graded $S$-module. Then, 
\begin{align}
H^i_{\M_l}(F)_{\n + \sign(1-i)\p} = 0 \label{2lemmaeq0}
\end{align}
for any $1 \le l \le k$, $i \ge 0$, $\p \in \kset{1-i}_k$ and $\n \in \ZZ^k$ such that $n_l \ge \rreg{l}(F)$.
\end{lemma}

\begin{proof} Let $\psi: \ZZ^k \rightarrow \ZZ^2$ be the group homomorphism given by $\n \mapsto (n_l, \sum_{i \not= l} n_i)$. Then $S^\psi$ is a standard $\NN^2$-graded ring and $F^\psi$ is a finitely generated $\ZZ^2$-graded $S^\psi$-module.
Note that under this $\NN^2$-graded structure of $S^\psi$, $\M_l$ is the $S^\psi$-ideal generated by $S^\psi_{(1,0)}$. Let
\begin{align}
\GG: 0 \rightarrow G_s \rightarrow \dots \rightarrow G_0 \rightarrow F \rightarrow 0 \label{2lemmares}
\end{align}
be a minimal $\ZZ^k$-graded free resolution of $F$. Since for each $i$, $G_i^\psi$ is a $\ZZ^2$-graded free $S^\psi$-module, $\GG^\psi$ is also a minimal $\ZZ^2$-graded free resolution of $F^\psi$ over $S^\psi$. For simplicity, let $G_i^\psi = \bigoplus_j S^\psi(-c_{ij}, -d_{ij})$ for $i = 0, \dots, s$. Let $(r_1, r_2) = \Rreg(F^\psi)$. Then, it follows from the definition of $\psi$ that $r_1 = \rreg{l}(F)$. For each $i=0, \dots, s$, let $c_i = \max_j \{ c_{ij} \}$. By Remark \ref{resolution-def}, we have
\begin{align}
r_1 = \max_i \{ c_i - i \}. \label{2lemmaduality}
\end{align}
Break (\ref{2lemmares}) into short exact sequences of the form
\begin{align}
0 \rightarrow M_{t+1} \rightarrow G_t^\psi \rightarrow M_t \rightarrow 0. \label{2lemmaseq}
\end{align}
Note that $M_0 = F^\psi$ and $M_s = G_s^\psi$. (\ref{2lemmaeq0}) is accomplished by showing that
\begin{align}
H^i_{\M_l}(M_t)_{(m-i+1+t, n)} & = 0 \ \forall \ 0 \le t \le s, i \ge 0, n \in \ZZ, m \ge r_1 \label{2lemmaind}
\end{align}
(when $t = 0$ this gives (\ref{2lemmaeq0})). We shall use descending induction on $t$ to prove (\ref{2lemmaind}).

For $t = s$, since local cohomology commutes with direct sums, we only need to show that $H^i_{\M_l}\big( S^\psi(-c_{sj}, -d_{sj}) \big)_{(m-i+1+s, n)} = 0$ for all $j$ and $i \ge 0, n \in \ZZ, m \ge r_1$, i.e. $$H^i_{\M_l}(S^\psi)_{(m-i+1+s-c_{sj}, n-d_{sj})} = 0 \ \forall \ j, i \ge 0, n \in \ZZ \ \text{and} \ m \ge r_1.$$
Let $\lambda: \ZZ^2 \rightarrow \ZZ$ be the group homomorphism given by $(n_1,n_2) \mapsto n_2$. Observe that under $\lambda$, $(S^\psi)^{[2]}_0 = \k[S^\psi_{(1,0)}]$ and $(S^\psi)^{[2]}_{n-d_{sj}}$ is a $\ZZ$-graded free $\k[S^\psi_{(1,0)}]$-module. Let $\mm_l$ be the maximal homogeneous ideal of $\k[S^\psi_{(1,0)}]$ (i.e. $\mm_l = \M_l \cap \k[S^\psi_{(1,0)}]$). Then, by Lemma \ref{2lem00}, we have
$$H^i_{\M_l}(S^\psi)_{(m-i+1+s-c_{sj}, n-d_{sj})} = H^i_{\mm_l}\big( (S^\psi)^{[2]}_{n-d_{sj}} \big)_{m-i+1+s-c_{sj}}.$$
It is easy to see (cf. \cite[Proposition 4.3]{hh}) that as a $\ZZ$-graded module over the polynomial ring $\k[S^\psi_{(1,0)}]$, $(S^\psi)^{[2]}_q$ is $0$-regular for any $q \in \ZZ$ (since it is a direct sum of copies of $\k[S^\psi_{(1,0)}]$). Thus, $H^i_{\mm_l}\big( (S^\psi)^{[2]}_{n-d_{sj}} \big)_{u-i+1} = 0$ for any $i \ge 0$ and $u \ge 0$. Moreover, by (\ref{2lemmaduality}), we have $r_1+s - c_{sj} \ge (c_s-s)+s-c_{sj} \ge 0$ for all $j$. Hence, with $m \ge r_1$, we have $H^i_{\mm_l}\big( (S^\psi)^{[2]}_{n-d_{sj}} \big)_{m-i+1+s-c_{sj}} = 0$ for any $i \ge 0$, i.e. $H^i_{\M_l}(S^\psi)_{(m-i+1+s-c_{sj}, n-d_{sj})} = 0 \ \forall \ i \ge 0$. This proves (\ref{2lemmaind}) for $t = s$.

Suppose that (\ref{2lemmaind}) has been proved for $t+1$. We shall show it for $t$. Taking the long exact sequence of local cohomology modules from the short exact sequence (\ref{2lemmaseq}) and since local cohomology commutes with direct sum, we have
$$\bigoplus_j H^i_{\M_l}\big( S^\psi(-c_{tj}, -d_{tj}) \big)_{(m-i+1+t,n)} \rightarrow H^i_{\M_l}(M_t)_{(m-i+1+t,n)} \rightarrow H^{i+1}_{\M_l}(M_{t+1})_{(m-i+1+t,n)}$$
By induction hypothesis, $H^{i+1}_{\M_l}(M_{t+1})_{(m-i+1+t,n)} = 0$ for any $i \ge 0, n \in \ZZ$ and $m \ge r_1$. By exactly the same line of arguments as in the case where $t = s$, we can show that $H^i_{\M_l}\big( S^\psi(-c_{tj}, -d_{tj}) \big)_{(m-i+1+t,n)} = 0$ for any $j, i \ge 0, n \in \ZZ$ and $m \ge r_1$. Thus, the above long exact sequence of local cohomology modules gives
$$H^i_{\M_l}(M_t)_{(m-i+1+t,n)} = 0 \ \forall \ i \ge 0, n \in \ZZ \ \text{and} \ m \ge r_1.$$
(\ref{2lemmaind}) is true by induction. Hence, (\ref{2lemmaeq0}) is proved.
\end{proof}

The following lemma is similar to Lemma \ref{a-inv-lem1}, but the use of Lemma \ref{hyrylemma} in the proof is no longer applicable because of the shifts $\sign(1-i)\p$ in the local cohomology modules. We shall need a new approach to prove the statement.

\begin{lemma} \label{2lem1}
Let $S$ be a standard $\NN^k$-graded polynomial ring over a field $\k$ and let $F$ be a finitely generated $\ZZ^k$-graded $S$-module. Let $\delta = \operatorname{proj-dim} F$. Then, for any nonempty index set $I \subseteq \{ 1, \dots, k \}$, we have
\begin{align}
H^i_{\M_I}(F)_{\n - \q + \sign(1-i)\p} = 0 \label{2lemmaeq1}
\end{align}
for $i \ge 0$, $\p \in \kset{1-i}_k, \q \in \kset{\delta}_k$ and $\n \in \ZZ^k$ such that $|\n_I| \ge \sum_{l \in I} \rreg{l}(F)+\delta|I|$.
\end{lemma}

\begin{proof} Let $r = |I|$. If $r = 1$ then (\ref{2lemmaeq1}) follows from Lemma \ref{2lem0}. Suppose $r \ge 2$. Let $\GG: 0 \rightarrow G_\delta \rightarrow \dots \rightarrow G_0 \rightarrow F \rightarrow 0$
be a minimal $\ZZ^k$-graded free resolution of $F$ over $S$, where $G_i = \bigoplus_j S(-\c{1}_{ij}, \dots, -\c{k}_{ij})$. For each $l = 1, \dots, k$ and each $i = 0, \dots, \delta$, set $\c{l}_i = \max_j \c{l}_{ij}$.

Consider first the case where $r = k$. In this case, $S^\phi$ is a standard $\NN$-graded ring with the maximal homogeneous ideal $\M_I = \M$ and $F^\phi$ is a finitely generated $\ZZ$-graded $S^\phi$-module. Also, $H^i_{\M}(F)^\phi$ is a $\ZZ$-graded $S$-module and $\GG^\phi$ is a minimal $\ZZ$-graded free resolution of $F^\phi$. It is easy to see that the highest shift in $G_i^\phi$ is at most $\sum_{l = 1}^k \c{l}_i$ for any $i = 0, \dots, \delta$. Thus,
$\reg(F^\phi) \le \max_{0 \le i \le \delta} \big\{ \sum_{l = 1}^k \c{l}_i - i \big\}$
where $\reg(F^\phi)$ is the usual $\ZZ$-graded Castelnuovo-Mumford regularity of $F^\phi$. Let $\n \in \ZZ^k$ be such that $|\n| \ge \sum_{l=1}^k \rreg{l}(F) + \delta k$. Then for any $0 \le i \le \delta$, we have $|\n|-\delta \ge \sum_{l=1}^k \rreg{l}(F) + \delta k - \delta \ge \sum_{l=1}^k (\c{l}_i - i) + \delta(k-1) = \sum_{l=1}^k \c{l}_i + (\delta-i)(k-1) - i \ge \sum_{l=1}^k \c{l}_i - i$. Thus, $|\n|-\delta \ge \reg(F^\phi)$. (\ref{2lemmaeq1}) now follows from the decomposition $[H^i_{\M}(F)^\phi]_m = \bigoplus_{|\m| = m} H^i_{\M}(F)_\m$.

Suppose $r < k$. Let $\psi: \ZZ^k \rightarrow \ZZ^2$ be the group homomorphism given by $\n \mapsto (|\n_I|, |\n_{\widehat{I}}|)$. Then $S^\psi$ is a standard $\NN^2$-graded ring, $F^\psi$ is a finitely generated $\ZZ^2$-graded $S^\psi$-module, $H^i_\M(F)^\psi$ and $H^i_{\M_I}(F)^\psi$ are $\ZZ^2$-graded $S^\psi$-modules, and $\GG^\psi$ is a minimal $\ZZ^2$-graded resolution of $F^\psi$. Let $(r_1, r_2)$ be the resolution regularity vector of $F^\psi$. Under the $\NN^2$-graded structure of $S^\psi$, $\M_I$ is the ideal generated by $S^\psi_{(1,0)}$.

Fix $\q \in \kset{\delta}_k$ and let $p_0 = \sum_{l \in I} \rreg{l}(F)+\delta|I|- |\q_I|$. Similar to the case when $r = k$, observe that $r_1 \le \max_{0 \le i \le \delta} \big\{ \sum_{l \in I} \c{l}_i - i \big\}$. Furthermore, for any $0 \le i \le \delta$, we have $p_0 \ge \sum_{l \in I} (\c{l}_i-i) + \delta|I|- |\q_I| \ge \sum_{l \in I} (\c{l}_i-i) + \delta|I| - \delta = \sum_{l \in I} \c{l}_i + (\delta-i)(|I|-1) - i \ge \sum_{l \in I} \c{l}_i - i$. Therefore, $p_0 \ge r_1$. It now follows from Lemma \ref{2lem0} that
$[H^i_{\M_I}(F)^\psi]_{(p-i+1,q)} = 0 \ \forall \ i \ge 0, q \in \ZZ \ \text{and} \ p \ge p_0.$
This implies that
$H^i_{\M_I}(F)_{\m} = 0 \ \forall \ i \ge 0, |\m_I| = p-i+1 \ \text{and} \ p \ge p_0.$
In particular, by taking $\m = \n - \q$, we have
$H^i_{\M_I}(F)_{\n - \q + \sign(1-i)\p} = 0 \ \text{for any} \ i \ge 0, \p \in \kset{1-i}_k$ and $\n \in \ZZ^k$ such that $|\n_I| \ge \sum_{l \in I} \rreg{l}(F)+\delta|I|$. This is true for any $\q \in \kset{\delta}_k$, so
(\ref{2lemmaeq1}) is proved.
\end{proof}

The following result, similar to Theorem \ref{a-inv-lemma}, gives the vanishing of multigraded pieces of local cohomology with respect to different multigraded ideals.

\begin{theorem} \label{2lemma}
Let $S$ be a standard $\NN^k$-graded polynomial ring over a field $\k$ and let $F$ be a finitely generated $\ZZ^k$-graded $S$-module. Let $\delta = \operatorname{proj-dim} F$. Then, for any $\n \ge \Rreg(F)+\delta\1, \q \in \kset{\delta}_k, i \ge 0, \p \in \kset{1-i}_k$ and any two disjoint index sets $I, J \subseteq \{ 1, \dots, k \}$, we have
\begin{align}
H^i_{\M_{I,J}}(F)_{\n - \q + \sign(1-i)\p} = 0. \label{2lemmaeq}
\end{align}
\end{theorem}

\begin{proof} The proof follows the same line of argument as that of Theorem \ref{a-inv-lemma} with a slight modification at the end. We shall use induction on $r = |I \cup J|$, the total number of elements in $I$ and $J$. If $r = 0$, then the statement is vacuous. If $r = 1$, then the statement clearly follows from Lemma \ref{2lem0}. Suppose $r \ge 2$. Without loss of generality, assume that $I = \{1, \dots, s\}$ and $J = \{s+1, \dots, r\}$. If $s = 0$, then the statement follows from Lemma \ref{2lem1}. Suppose $s > 0$. Let $K = I \backslash \{s\}$ and $L = J \cup \{s\}$. By identifying $\M_{I, \emptyset}$ with $\M_{K, \{s\}}$, we may assume that $J \not= \emptyset$, i.e. $s < r$. We now use ascending induction on $s$. Observe that $\M_{K,J} + \M_{I, \emptyset} = \M_{K,L}$ and $\M_{K,J} \cap \M_{I, \emptyset} = \M_{I,J}$. Consider the following Mayer-Vietoris sequence of local cohomology
$$\dots \rightarrow H^i_{\M_{K,J}}(F) \oplus H^i_{\M_{I, \emptyset}}(F) \rightarrow H^i_{\M_{I,J}}(F) \rightarrow H^{i+1}_{\M_{K,L}}(F) \rightarrow \dots$$
By induction on $s$, we have
$H^{i+1}_{\M_{K,L}}(F)_{\n-\q+\sign(-i)\p'} = 0$
for all $\n \ge \Rreg(F)+\delta\1, \q \in \kset{\delta}_k, i \ge 0$ and $\p' \in \kset{-i}_k$. 

For any $\p \in \kset{1-i}_k$, if $i = 0$ then take $\p' = \0$ and if $i > 0$ then take $\p' = \p+\e_j$ for some $j$. Then, $\p' \in \kset{-i}_k$. Moreover, if $i=0$ then $\n+\sign(1-i)\p = (\n+\p) + \sign(-i)\p'$ and if $i > 0$ then $\n+\sign(1-i)\p = (\n+\e_j) + \sign(-i)\p'$. Thus, we have $H^{i+1}_{\M_{K,L}}(F)_{\n - \q + \sign(1-i)\p} = 0$ for all $\n \ge \Rreg(F)+\delta\1, \q \in \kset{\delta}_k, i \ge 0$ and $\p \in \kset{1-i}_k$. By induction on $r$, we also have
$H^i_{\M_{K,J}}(F)_{\n - \q + \sign(1-i)\p} = 0 = H^i_{\M_{I, \emptyset}}(F)_{\n - \q + \sign(1-i)\p}$
for all $\n \ge \Rreg(F)+\delta\1, \q \in \kset{\delta}_k, i \ge 0$ and $\p \in \kset{1-i}_k$. Hence,
$H^i_{\M_{I,J}}(F)_{\n - \q + \sign(1-i)\p} = 0$
for all $\n \ge \Rreg(F)+\delta\1, \q \in \kset{\delta}_k, i \ge 0$ and $\p \in \kset{1-i}_k$. (\ref{2lemmaeq}) is proved.
\end{proof}


\subsection{Multigraded regularity and the resolution regularity vector.} \label{reg-subsec} We proceed to prove our next main theorem, which provides a lower and an upper bound for multigraded regularity in terms of the resolution regularity vector and other $\ZZ$-graded invariants. We also give an example showing that the inclusion $\Rreg(F)+\NN^k \subseteq \Reg(F)$ does not always hold in general.

\begin{theorem} \label{2regbound}
Let $S$ be a standard $\NN^k$-graded polynomial ring over a field $\k$ and let $F$ be a finitely generated $\ZZ^k$-graded $S$-module. Let $\delta = \operatorname{proj-dim} F$. For each $l = 1, \dots, k$, let $\lr{l}$ be the least integer such that
$\la{l}^{\lr{l}}(F) = \la{l}^*(F).$
Let $\r(F) = \big(\la{1}^*(F)+\lr{1}, \dots, \la{k}^*(F)+\lr{k}\big)$. Then,
\begin{align}
\bigcup_{\q \in \kset{\delta}_k} \big( \Rreg(F) + \delta\1 - \q + \NN^k \big) \subseteq \Reg(F) \subseteq \r(F) + \NN^k. \label{2regboundeq}
\end{align}
\end{theorem}

\begin{proof} Pick $I = \{ 1, \dots, k \}$ and $J = \emptyset$. Then, $\M_{I,J} = \B$. It follows from Theorem \ref{2lemma} that for any $\n \ge \Rreg(F)+\delta\1, i \ge 0$, $\p \in \kset{1-i}_k$ and $\q \in \kset{\delta}_k$, we have
$$H^i_{\B}(F)_{\n - \q + \sign(1-i)\p} = 0.$$
By definition, this implies that $\Rreg(F)+\delta\1 -\q \in \Reg(F)$ for any $\q \in \kset{\delta}_k$. The left hand inclusion of (\ref{2regboundeq}) is proved.

Fix $\n \in \Reg(F)$ and $l \in \{ 1, \dots, k \}$. To prove the right hand inclusion of (\ref{2regboundeq}), we need to show that $n_l \ge \la{l}^*(F)+\lr{l}$. 
Let $\m \in \ZZ^k$ be such that $\m_{\widehat{l}} \gg \0$ and $m_l = n_l$. In particular, we shall pick $\m \ge \n$ so that $\m \in \Reg(F)$.

Consider first the case when $\lr{l} \ge 2$. It follows from the proof of Proposition \ref{2projection} that $\m^{[l]} \not\in \aa^{\lr{l}}(F)+\1$. This implies that $\m^{[l]} + (\lr{l}-1)\e_l \not\in \aa^{\lr{l}}(F)+\1+(\lr{l}-1)\e_l$. Moreover $(\aa^{\lr{l}}(F)+\1+(\lr{l}-1)\e_l)$ is a term in the intersection $\bigcap_i \bigcap_{\p \in \kset{1-i}_k} \big( \aa^i(F)+\1 - \sign(1-i)\p\big).$ Thus, $\m^{[l]} + (\lr{l}-1)\e_l \not\in \bigcap_i \bigcap_{\p \in \kset{1-i}_k} \big( \aa^i(F)+\1 - \sign(1-i)\p\big)$. By Proposition \ref{2a-invariant}, $\m^{[l]}+(\lr{l}-1)\e_l \not\in \Reg(F)$. Now, if $m_l \le \la{l}^*(F)+\lr{l}-1$ then $\m^{[l]}+(\lr{l}-1)\e_l \ge \m$, and so, $\m^{[l]}+(\lr{l}-1)\e_l \in \Reg(F)$ which is a contradiction. Thus, we must have $n_l = m_l \ge \la{l}^*(F)+\lr{l}$.
On the other hand, suppose $\lr{l} \le 1$. If $m_l < \la{l}^*(F)+\lr{l} \le \la{l}^*(F)+1$ then $\m^{[l]} \ge \m$, and thus, $\m^{[l]} \in \Reg(F)$. By Remark \ref{regvsainv}, we have $\m^{[l]} \in (\aa^0(F) + \1) \bigcap (\aa^1(F)+\1)$. However, it follows from the proof of Proposition \ref{2projection} that either $\m^{[l]} \not\in \aa^0(F)+\1$ or $\m^{[l]} \not\in \aa^1(F)+\1$, a contradiction. Therefore, we must also have $n_l = m_l \ge \la{l}^*(F)+\lr{l}$.
\end{proof}

\begin{example} \label{example5}
It follows from Example \ref{example2} and the minimal free resolution $0 \rightarrow S \rightarrow S \rightarrow 0$ of $S$ over itself that $\Reg(S) = \Rreg(S) + \NN^k = \NN^k$. This shows that the left hand inclusion of (\ref{2regboundeq}) is sharp. We shall see later, in Example \ref{example6}, that the right hand inclusion of (\ref{2regboundeq}) is also sharp.
\end{example}

The following example gives a negative answer to \cite[Question 2.6]{sv}.

\begin{example} \label{resnotinreg}
Let $S = \k[x_1, \dots, x_N, y_1, \dots, y_N]$ be the coordinate ring of $\PP^{N-1} \times \PP^{N-1}$. Let $f_1, \dots, f_N$ be a regular sequence in $S$ such that $\operatorname{bi-degree}(f_i) = (1,1)$ (we can construct $f_1, \dots, f_N$ as follows: take $g_1, \dots, g_N \in \k[x_1, \dots, x_N]$ and $h_1, \dots, h_N \in \k[y_1, \dots, y_N]$ to be regular sequences, and take $f_i = g_ih_i$ for $i = 1, \dots, N$). Let $F = S/(f_1, \dots, f_N)$. It can be observed that the Koszul complex of $f_1, \dots, f_N$ gives a minimal $\ZZ^2$-graded free resolution of $F$ over $S$. It then follows that $\Rreg(F) = \0$ (since the Koszul complex of $f_1, \dots, f_N$ is linear in each coordinate of the degree vector). By Lemma \ref{2lem0}, we have
\begin{align}
H^i_{\M_1}(F)_{\n+\sign(1-i)\p} = H^i_{\M_2}(F)_{\n+\sign(1-i)\p} = 0 \label{examsummand}
\end{align}
for any $\n \ge \0$ and $\p \in \kset{1-i}_k$.

As before, let $\phi: \ZZ^2 \rightarrow \ZZ$ be the group homomorphism given by $(m_1,m_2) \mapsto m_1+m_2$. Then $F^\phi$ is a $\ZZ$-graded $S^\phi$-module. It follows from the Koszul complex of $f_1, \dots, f_N$ (since $\deg(f_i) = 2 \ \forall \ i$) that $\reg(F^\phi) = N$. By (\ref{reg-ainv}), we have $\max_i \{ a_i(F^\phi) + i \} = N$. Moreover, since $f_1, \dots, f_N$ form a regular sequence, $F$ is a Cohen-Macaulay $S$-module of dimension $N$. This implies that $a_i(F^\phi) = -\infty$ for $i < N$. Thus, we have $a_N(F^\phi)+N = N$, i.e. $a_N(F^\phi) = 0$. By definition, this says that
\begin{align}
\bigoplus_{m_1+m_2 = 0}H^N_\M(F)_{(m_1,m_2)} = [H^N_\M(F)^\phi]_0 \not= 0. \label{exampleh}
\end{align}

On the other hand, by successively considering the exact sequence
$$0 \rightarrow S/(f_1, \dots, f_{i-1}) \stackrel{\times f_i}{\rightarrow} S/(f_1, \dots, f_{i-1}) \rightarrow S/(f_1, \dots, f_i) \rightarrow 0$$
and its long exact sequence of local cohomology modules with respect to $\M$, for $i = 1, \dots, N$, we deduce that $a_{l,*}(F) \le 0$ for $l = 1,2$. This implies that $\a^*(F) \le \0.$ That is, $H^i_\M(F)_{(m_1,m_2)} = 0$ if $m_1 > 0$ or $m_2 > 0$.
This, together with (\ref{exampleh}), implies that 
\begin{align}
H^N_\M(F)_{(0,0)} \not= 0. \label{examplehh}
\end{align}

Observe further that $\M = \M_1 + \M_2$ and $\B = \M_1 \cap \M_2$. Consider the Mayer-Vietoris sequence of local cohomology modules
\begin{align}
\dots \rightarrow H^i_\M(F) \rightarrow H^i_{\M_1}(F) \oplus H^i_{\M_2}(F) \rightarrow H^i_{\B}(F) \rightarrow \dots \label{mvsequence}
\end{align}
It follows from (\ref{examsummand}) and (\ref{examplehh}) that $H^N_{\B}(F)_{(0,0)} \not= 0$ (by choosing $\n = \p \in \kset{1-N}_2$ in (\ref{examsummand}) and applying (\ref{examplehh})). This implies that $(a,b) \not\in \Reg(F)$ for any $(a,b) \in \NN^2$ such that $a+b = N-1$. This shows, in a way, that the left hand inclusion of (\ref{2regboundeq}) is sharp. In particular, we have $(N-1, 0) \not\in \Reg(F)$ and $(0,N-1) \not\in \Reg(F)$. More interestingly, it also follows that $\Rreg(F) = (0,0) \not\in \Reg(F)$ (since $(N-1,0) > (0,0)$). This shows that the inclusion $\Rreg(F) + \NN^k \subseteq \Reg(F)$ does not always hold in general.
\end{example}

The following corollary will be useful in the next section.

\begin{corollary} \label{regvslocalainv}
Using the same notation as in Theorem \ref{2regbound}. Then, for each $l = 1, \dots, k$, we have
$$\rreg{l}(F) > \la{l}^*(F).$$
\end{corollary}

\begin{proof} By taking $\q = \delta\e_l$, Theorem \ref{2regbound} implies that $\Rreg(F)+\delta(\1-\e_l) \in \Reg(F)$. This, together with Remark \ref{regvsainv} and Theorem \ref{a-inv-bound}, implies that $\Rreg(F)+\delta(\1-\e_l) \in \aa^*(F)+\1 \subseteq \La^*(F)+\1+\NN^k$. Thus, $\Rreg(F)+\delta(\1-\e_l) \ge \La^*(F)+\1$. Hence, $\rreg{l}(F) \ge \la{l}^*(F)+1 > \la{l}^*(F)$. The statement is proved.
\end{proof}

We shall now discuss a number of important special cases of Theorem \ref{2regbound}. Suppose $F$ is a $\ZZ^k$-graded $S$-module and $\F$ its associated coherent sheaf on $X$. We define the {\it support} of $F$ to be
\begin{align*}
\suppo(F) = \big\{ \wp \in X \ \big| \ F_{(\wp)} \not= 0 \big\} = \big\{ \wp \in X \ \big| \ \F_\wp \not= 0 \big\}
\end{align*}
where $F_{(\wp)}$ is the $\ZZ^k$-graded homogeneous localization of $F$ at $\wp$, and $\F_\wp$ denotes the stalk of $\F$ at $\wp$. For each $l = 1, \dots, k$, let $\supp{l}(F) = \ppi{l}\big(\suppo(F)\big)$.

Inspired by the $\ZZ$-graded situation, we define the {\it fiber regularity vector} of $F$ as follows.

\begin{defn} \label{localregvector} For each $l = 1, \dots, k$, let
$$\lreg{l}(F) = \max_i \big\{ \la{l}^i(F) + i \big\}.$$
The {\it fiber regularity vector} of $F$ is defined to be the vector
$$\Lreg(F) = (\lreg{1}(F), \dots, \lreg{k}(F)) \in \ZZ^k.$$
\end{defn}

When $F$ is a Cohen-Macaulay $S$-module we obtain a particularly nice upper bound on multigraded regularity in terms of the fiber regularity vector.

\begin{corollary} \label{bound-cor}
Let $S$ be a standard $\NN^k$-graded polynomial ring over a field $\k$ and let $F$ be a finitely generated Cohen-Macaulay $\ZZ^k$-graded $S$-module. Then,
$$\Reg(F) \subseteq \Lreg(F) + \NN^k.$$
\end{corollary}

\begin{proof} We first observe that for any two $\NN^{k-1}$-graded homogeneous prime ideals $\qq \subseteq \pp \subseteq \mS{l}$, we have $H^i_{S_{(\qq)+}}(F_{(\qq)}) = \big[H^i_{S_{(\pp)+}}(F_{(\pp)})\big]_{(\qq)}$. Thus, $a_i(F_{(\qq)}) \le a_i(F_{(\pp)})$ for any $i \ge 0$. This implies that $\la{l}^i(F)$ is in fact given by
$$\la{l}^i(F) = \max \big\{ a_i(F_{(\pp)}) \ \big| \ \pp \in \X{l} \ \text{a closed point} \big\} \ \forall \ l = 1, \dots, k.$$
It can further be seen that since $F$ is a Cohen-Macaulay $S$-module, $F_{(\pp)}$ is a Cohen-Macaulay $S_{(\pp)}$-module for any $l = 1, \dots, k$ and $\pp \in \X{l}$. That is, $\la{l}^i(F) = 0$ for any $i < \dim \supp{l}(F)$ (since for a closed point $\pp \in \X{l}$, $\dim F_{(\pp)} = \dim \supp{l}(F)$). Clearly, $\la{l}^i(F) = 0$ for $i > \dim \supp{l}(F)$. By definition, we now have $\la{l}^*(F) = \la{l}^{\dim \supp{l}(F)}(F)$. Using the same notations as in Theorem \ref{2regbound}, let $\lr{l}$ be the least integer such that $\la{l}^*(F) = \la{l}^{\lr{l}}(F)$, then we have $\lr{l} = \dim \supp{l}(F)$. Hence, $\lreg{l}(F) = \la{l}^{\dim \supp{l}(F)}(F) + \dim \supp{l}(F) = \la{l}^*(F) + \lr{l}$. The conclusion now follows from that of Theorem \ref{2regbound}.
\end{proof}

\begin{example} \label{example6}
As in Example \ref{example4}, we observe that for each $l \in \{ 1, \dots, k \}$ and $\pp \in \X{l}$, we have $\reg\big(S_{(\pp)}\big) = 0$. This implies that $\rreg{l}(S) \ge 0$ for any $l = 1, \dots, k$. Thus, $\Lreg(F) \ge \0$. Moreover, by Corollary \ref{bound-cor} and Example \ref{example2}, we have $\NN^k \subseteq \Lreg(F) + \NN^k$, i.e. $\Lreg(F) \le \0$. Hence, $\Lreg(F) = \0$. This shows that the right hand inclusion of Corollary \ref{bound-cor} is sharp.
\end{example}

\begin{question} Does the inclusion $\Reg(F) \subseteq \Lreg(F) + \NN^k$ always hold?
\end{question}

Next, we address a situation where $F$ has a good minimal $\ZZ^k$-graded free resolution. In this situation, we shall see that the resolution regularity vector of $F$ coincides with its local regularity vector.

\begin{defn} \label{gooddef} Let $S$ be a standard $\NN^k$-graded polynomial ring over a field $\k$ and let $F$ be a finitely generated $\ZZ^k$-graded $S$-module. Let $\GG: 0 \rightarrow G_s \rightarrow \dots \rightarrow G_0 \rightarrow F \rightarrow 0$ be a minimal $\ZZ^k$-graded free resolution of $F$ over $S$, where $G_i = \bigoplus_j S(-\n^{i,j})$. For each $0 \le i \le s$ and each $1 \le l \le k$, let $n^i_l = \max_j \{\n^{i,j}_l\}$ be the maximum of the $l$th coordinate of the shifts $\n^{i,j}$. We shall say that $F$ has {\it purely increasing minimal $\ZZ^k$-graded free resolutions} if $n^i_l > n^{i-1}_l$ for all $1 \le i \le s$ and $1 \le l \le k$.
\end{defn}
Note that since all minimal $\ZZ^k$-graded free resolutions of $F$ are isomorphic, our property of having purely increasing minimal $\ZZ^k$-graded free resolutions is well-defined and does not depend on a particular minimal free resolution $\GG$. When $F$ has purely increasing minimal $\ZZ^k$-graded free resolutions, we will also call any minimal $\ZZ^k$-graded free resolution of $F$ a purely increasing minimal $\ZZ^k$-graded free resolution.

The following result allows us to find $\Reg(F)$ by checking only a finite number of possibilities if it is known that $F$ possesses a good minimal $\ZZ^k$-graded free resolution.

\begin{corollary} \label{equality}
Let $S$ be a standard $\NN^k$-graded polynomial ring over a field $\k$ and let $F$ be a finitely generated Cohen-Macaulay $\ZZ^k$-graded $S$-module with purely increasing minimal $\ZZ^k$-graded free resolutions. Let $\delta = \operatorname{proj-dim} F$. Then, we have
$$\Reg(F) \subseteq \Rreg(F) + \NN^k.$$
\end{corollary}

\begin{proof} In view of Corollary \ref{bound-cor}, it suffices to show that
$\Rreg(F) = \Lreg(F).$ 

Let $\GG: 0 \rightarrow G_s \rightarrow \dots \rightarrow G_0 \rightarrow F \rightarrow 0$ be a minimal $\ZZ^k$-graded free resolution of $F$ over $S$, where $G_i = \bigoplus_j S(-\n^{i,j})$. Fix an integer $l \in \{1, \dots, k\}$. For any $\pp \in \X{l}$, it can be seen that the extension $S \hookrightarrow S_\pp$ is a flat extension. Thus, $\GG_\pp = \GG \otimes_{\mS{l}} \mS{l}_\pp$ is a $\ZZ^k$-graded free resolution of $F_\pp$ over $S_\pp$. By taking elements of degree $\0 \in \ZZ^{k-1}$ of $\GG_\pp^{\phi_{\widehat{l}}}$, we obtain a $\ZZ$-graded free resolution of $F_{(\pp)}$, namely 
$$\GG_{(\pp)}: 0 \rightarrow \big[ G_s \big]_{(\pp)} \rightarrow \dots \rightarrow \big[ G_0 \big]_{(\pp)} \rightarrow F_{(\pp)} \rightarrow 0.$$
Here, $\big[ G_i \big]_{(\pp)} = \bigoplus_j \big[ S^{\phi_{\widehat{l}}}(-\n^{i,j}_{\widehat{l}}) \big]_{(\pp)}\big(-\n^{i,j}_l\big) = \bigoplus_j S_{(\pp)}(-\n^{i,j}_l)$. 

Observe that $\GG_{(\pp)}$ is not necessarily minimal. In fact, in general, there will be cancelations among the direct summands of $\big[ G_i \big]_{(\pp)}$'s that would result in a minimal free resolution of $F_{(\pp)}$. However, we shall show that in these cancelations, the highest shifts in $\big[ G_i \big]_{(\pp)}$'s remain the same. Indeed, let $M_i$ be the presentation matrix of $G_i \rightarrow G_{i-1}$ for $i \ge 0$. Since $\GG$ is a minimal resolution of $F$ over $S$, the entries of $M_i$ are all in $\M$. As before, let $n^i_l = \max_j \{\n^{i,j}_l\}$ and suppose $j_0$ is an index that $\n^{i,j_0}_l = n^i_l$. Observe further that since $\GG$ is a purely increasing resolution of $F$, we have $n^i_l > n^{i-1}_l$. This implies that entries of the $j_0$-th column of $M_i$ are all in $\M_l$. Moreover, the localization $S_\pp = S \otimes_{\mS{l}} \mS{l}_\pp$ is the operation that inverts elements in $\mS{l} \backslash \pp$. Thus, in the matrix $M_i \otimes_{\mS{l}} \mS{l}_\pp$, entries of the $j_0$-th column are in $S_{(\pp)} \cap \big( M_i \otimes_{\mS{l}} \mS{l}_\pp \big) = S_{(\pp)+}$. This means that the shift $\n^{i,j_0}_l$ is not cancelled in the minimal free resolution of $F_{(\pp)}$. 

We have shown that the highest shift in $\big[ G_i \big]_{(\pp)}$ is exactly $\max_j\{\n^{i,j}_l\}$ for any $i \ge 0$. This together with Remark \ref{resolution-def} and its well known $\ZZ$-graded version imply that $\rreg{l}(F) = \operatorname{reg}(F_{(\pp)})$ for $\pp \in \X{l}$. Hence, $\rreg{l}(F) = \lreg{l}(F)$. Since this is true for any $l = 1, \dots, k$, we have $\Rreg(F) = \Lreg(F)$.
\end{proof}

\begin{corollary} \label{ci}
Let $Y$ be a complete intersection subscheme of $X$ cut out by a regular sequence $f_1, \dots, f_s \in S$, where $\operatorname{multideg}(f_i) = (d_{i1}, \dots, d_{ik}) \ge \1 \ \forall \ i$. Let $F$ be the $\NN^k$-graded coordinate ring of $Y$ and let $\delta = \operatorname{proj-dim} F$. For $l = 1, \dots, k$, let
$$\dd_l = \max_{0 \le i \le s} \big\{ \max \big\{ \sum_{j \in J} d_{jl} \ \big| \ J \subseteq \{ 1, \dots, k \}, |J| = i \big\} - i \big\}.$$
Let $\dd = (\dd_1, \dots, \dd_k)$. Then,
$\bigcup_{\q \in \kset{\delta}_k} \big(\dd + \delta\1 - \q + \NN^k\big) \subseteq \Reg(F) \subseteq \dd + \NN^k.$
\end{corollary}

\begin{proof} As before, we first observe that a minimal $\ZZ^k$-graded free resolution of $F$ is given by the Koszul complex associated to the regular sequence $f_1, \dots, f_s$. It follows from the condition $\operatorname{multideg}(f_i) \ge \1 \ \forall \ i$ that the Koszul complex is a purely increasing minimal $\ZZ^k$-graded resolution of $F$. Thus, $F$ has purely increasing minimal $\ZZ^k$-graded free resolutions. It also follows from the Koszul complex that $\Rreg(F) = \dd$. Moreover, since $f_1, \dots, f_s$ form a regular sequence, $F$ is a Cohen-Macaulay $S$-module. The statement now follows from Corollary \ref{equality}.
\end{proof}

\begin{remark} \label{strictnoteq}
It follows from Corollary \ref{ci} and our calculations in Example \ref{resnotinreg} that when $S$ is the coordinate ring of $\PP^{N-1} \times \PP^{N-1}$ and $F = S/(f_1, \dots, f_N)$, where $f_1, \dots, f_N$ is a regular sequence in $S$, we have a strict inclusion
$\Reg(F) \subsetneq \Rreg(F) + \NN^2.$
\end{remark}


\subsection{The resolution regularity vector and global generation.} \label{res-subsec} We conclude the paper by proving an important property of the resolution regularity vector, which shows that the resolution regularity shares similarities with the original definition of Castelnuovo-Mumford regularity.

\begin{lemma} \label{reslem1}
Let $S$ be a standard $\NN^k$-graded polynomial ring over a field $\k$ and let $F$ be a finitely generated $\ZZ^k$-graded $S$-module. Let $\F$ be the coherent sheaf associated to $F$ on $X = \proj S$. Then, for any $\m \ge \Rreg(F)$, we have
\begin{align}
H^0(X, \F(\m)) = F_\m. \label{reslem1eq}
\end{align}
\end{lemma}

\begin{proof} We shall prove (\ref{reslem1eq}) by induction in $k$. For $k = 1$, the statement is well known since $\Reg(F) = \Rreg(F) + \NN = \big\{ r \ \big| \ r \ge \reg(F) \big\}$ in this case. Suppose now that $k \ge 2$ and (\ref{reslem1eq}) has been proved for $k-1$. We shall show that (\ref{reslem1eq}) holds for $k$. Indeed, suppose $\m \ge \Rreg(F)$. By Corollary \ref{regvslocalainv}, we have $m_l > \la{l}^*(F)$. Thus, it follows from Proposition \ref{2projection} that $\ppi{l}_* \F(\m) = \widetilde{\lF{l}_{m_l}}(\m_{\widehat{l}})$. Therefore,
\begin{align}
H^0(X, \F(\m)) = H^0(\X{l}, \ppi{l}_* \F(\m)) = H^0\big(\X{l}, \widetilde{\lF{l}_{m_l}}(\m_{\widehat{l}})\big). \label{reslem1proj1}
\end{align}
Let $\GG: 0 \rightarrow G_s \rightarrow \dots \rightarrow G_0 \rightarrow F \rightarrow 0$ be a minimal $\ZZ^k$-graded free resolution of $F$ over $S$. By taking elements of degree $m_l\e_l$ from $\GG$, we obtain a $\ZZ^{k-1}$-graded free resolution of $\lF{l}_{m_l}$ over $\mS{l}$. This resolution is not necessarily minimal but by considering its $\ZZ^{k-1}$-graded shifts, it follows from Remark \ref{resolution-def} that $\Rreg\big( \lF{l}_{m_l}\big) \le \Rreg(F)_{\widehat{l}} \le \m_{\widehat{l}}$. By induction, we now have
\begin{align}
H^0\big(\X{l}, \widetilde{\lF{l}_{m_l}}(\m_{\widehat{l}})\big) = \big[ \lF{l}_{m_l} \big]_{\m_{\widehat{l}}} = F_\m. \label{reslem1proj2}
\end{align}
It follows from (\ref{reslem1proj1}) and (\ref{reslem1proj2}) that (\ref{reslem1eq}) is true for $k$. The lemma is proved.
\end{proof}

\begin{theorem} \label{globalgen}
Let $S$ be a standard $\NN^k$-graded polynomial ring over a field $\k$ and let $F$ be a finitely generated $\ZZ^k$-graded $S$-module. Let $\F$ be the coherent sheaf associated to $F$ on $X = \proj S$. Then, for $\m \ge \Rreg(F)$, we have
\begin{enumerate}
\item The natural map $H^0(X, \F(\m)) \otimes H^0(X, \O_X(\n)) \rightarrow H^0(X, \F(\m+\n))$ is surjective for any $\n \in \NN^k$.
\item $\F(\m)$ is generated by its global sections.
\end{enumerate}
\end{theorem}

\begin{proof} Since $\m \ge \Rreg(F)$, $\m$ bounds the multi-degrees of a minimal system of generators for $F$. This implies that $F_\m S_\n = F_{\m+\n}$ for $\m \ge \Rreg(F)$ and $\n \in \NN^k$. It also follows from Lemma \ref{reslem1} that $H^0(X, \F(\m)) = F_\m$ and $H^0(X, \F(\m+\n)) = F_{\m+\n}$ for $\m \ge \Rreg(F)$ and $\n \in \NN^k$. Thus, (1) holds since $H^0(X, \O_X(\n)) = S_\n$ for $\n \in \NN^k$. To prove (2), we observe again that $\m$ bounds the multi-degrees of a minimal system of generators for $F$. Therefore, elements of $F_\m$ generate $\bigoplus_{\q \ge \m} F_\q$. This and the fact that $\F$ coincides with the coherent sheaf associated to $\bigoplus_{\q \ge \m} F_\q$ establish (2).
\end{proof}



\end{document}